\documentclass[trsc,nonblindrev]{informs3} % current default for manuscript submission

\usepackage{fancyhdr}
\def\vphind{\textcolor{black}{\varphi^n_{i}}}
\def\vphjmd{\textcolor{black}{\varphi^m_{j}}}
\def\vphjnd{\textcolor{black}{\varphi^n_{j}}}

\def\vphjdd{\textcolor{black}{\varphi^{\sbullet[.5]}_{j}}}
\def\vphldd{\textcolor{black}{\varphi^{\sbullet[.5]}_{\ell}}}
\def\Qcont{\textcolor{black}{Q^n_{i,j}}}

\def\bM{{\overline M}}
\def\wM{{\widetilde M}}
\def\uw{{\underline w}}
\newcommand\sbullet[1][.5]{\mathbin{\vcenter{\hbox{\scalebox{#1}{$\bullet$}}}}}
\def\vfb{\varphi^{\sbullet[.5]}}

\def\locf#1{\Bigl\lfloor#1 \Bigr\rfloor}

\def\one{\mathbf{1}}

\OneAndAHalfSpacedXI % current default line spacing
\usepackage{natbib}
 \bibpunct[, ]{(}{)}{,}{a}{}{,}%
 \def\bibfont{\small}%

\usepackage{amsmath, amssymb, bbm, amsfonts}

\usepackage{enumitem}
\usepackage[utf8]{inputenc}
\usepackage[UKenglish]{isodate}
\usepackage[font=small,labelfont=it]{caption}
\usepackage{subcaption}
\usepackage{multirow}
\usepackage{fancyvrb}
\usepackage[figuresleft]{rotating}
\usepackage{geometry}
\usepackage{hyperref}
\usepackage{booktabs}
\usepackage{makecell}

\usepackage[dvipsnames,table]{xcolor}
\definecolor{tableShade}{gray}{0.9}

\usepackage{graphicx}
\graphicspath{{img/}}
\TheoremsNumberedThrough     % Preferred (Theorem 1, Lemma 1, Theorem 2)
\EquationsNumberedThrough    % Default: (1), (2), ...
\begin{document}
%\setcounter{page}{1}
%\DOI{10.1287/xxxx.0000.0000}%

\RUNAUTHOR{C.-A. Lehalle and G. Livieri}

\RUNTITLE{Competition among seaports through Mean Field Games}

\TITLE{Competition among seaports through Mean Field Games and real-world data}
%\TITLE{}

\ARTICLEAUTHORS{%
\AUTHOR{Charles-Albert Lehalle}
\AFF{Polytechnique, CMAP, IP Paris and  Fellow of Institut Louis Bachelier, Paris, \EMAIL{charles-albert.lehalle@polytechnique.edu}}
\AUTHOR{Giulia Livieri}
\AFF{London School of Economics and Political Science, \EMAIL{g.livieri@lse.ac.uk}}
} % end of the block

\ABSTRACT{This paper presents a Mean Field Game (MFG) model for maritime traffic flow, treating the navigation of ships between seaports as a large-scale stochastic control problem. The MFG framework enables the modeling of agents at a microscopic level as rational decision-makers who seek to optimize their utility, thereby translating complex microscopic behaviors into macroscopic models. We build upon this MFG framework to develop a mesoscopic-scale MFG model that defines the payoff and cost functions for a coordinator at each seaport considered in our study. The coordinator determines the routes taken by ships transporting goods between ports by evaluating several key factors: transportation costs, expected profit margins from loading specific goods at the seaports and unloading them at various destinations, and a congestion term that reflects the costs associated with accessing the destination port. We derive an explicit solution for the stationary version of the model under certain approximations and establish conditions necessary to ensure the uniqueness of the corresponding Mean Field Equilibrium (MFE). Furthermore, we introduce a statistical methodology to infer the parameters of the game from real-world data, specifically focusing on costs and the components of expected commercial margins. To validate our model in a real-world context, we analyze the \texttt{ShipFix} dataset of daily ``Dry Coal" shipments worldwide from 2015 to 2025. Our discussion highlights the influence of empirical traffic flow on various components of costs. We believe that this research represents a significant advancement in the application of MFGs for effective maritime traffic management and offers valuable insights for practitioners in the field.}

%\KEYWORDS{deterministic inventory theory; infinite linear programming duality; 
%  existence of optimal policies; semi-Markov decision process; cyclic schedule}

\KEYWORDS{Mean Field Games, maritime traffic flow, congestion costs, network flows, parameters Inference, \texttt{ShipFix} dataset.}

\maketitle
%%%%%%%%%%%%%%%%%%%%%%%%%%%%%%%%%%%%%%%%%%%%%%%%%%%%%%%%%%%%%%%%%%%%%%

% Samples of sectioning (and labeling) in TRSC
% NOTE: (1) \section and \subsection do NOT end with a period
%       (2) \subsubsection and lower need end punctuation
%       (3) capitalization is as shown (title style).
%
%\section{Introduction.}\label{intro} %%1.
%\subsection{Duality and the Classical EOQ Problem.}\label{class-EOQ} %% 1.1.
%\subsection{Outline.}\label{outline1} %% 1.2.
%\subsubsection{Cyclic Schedules for the General Deterministic SMDP.}
%  \label{cyclic-schedules} %% 1.2.1
%\section{Problem Description.}\label{problemdescription} %% 2.

% Text of your paper here
\section{Introduction}\label{section::introduction}
Competition among seaports plays a crucial role in driving performance improvements within the shipping industry. It has become a significant issue in transportation economics, especially as the industry adapts to structural changes in the maritime supply chain. Maritime transportation is a vital component of global trade, with approximately $80\%$ of commerce by volume and $70\%$ by value being transported by sea and processed at seaports worldwide (e.g., \cite{unctad2023}). Maritime flow models that incorporate decision-making sources are essential tools for effectively evaluating the significant impacts of major global events. For instance, the effects of pandemics, such as COVID-19, and geopolitical conflicts, such as the disruption of goods at key seaports due to the war in Ukraine, highlight how quickly stability can be compromised (e.g., \cite{komaromi2022supply}). This includes the unavailability of essential commodities like wheat and fertilizers, which are primarily produced in Ukraine, as well as the direct closure of seaports. Moreover, extreme weather events, such as wildfires and hurricanes, can severely disrupt traffic across various regions. In this paper, we aim to develop a game-theoretic maritime flow model using Mean Field Games (MFGs). This model will enable us to establish the payoff and cost functions for a coordinator at each of the seaports considered. The coordinator will decide the routes taken by ships transporting goods between seaports on a mesoscopic scale.
\subsection{Motivation}\label{subsection::motivation}
Mean Field Game (MFG) is a game-theoretic framework that has gained popularity for designing decision-making processes in systems with many agents. It was introduced almost simultaneously by \cite{huang2006large} and \cite{lasry2007mean}. Most analyses of MFG have concentrated on continuous state and time symmetric differential games. In the MFG framework, each player seeks to solve an optimal control problem based on an objective function, while other players make similar decisions concurrently. Players interact on a microscopic level through their coupled objective functions. As the number of players increases, the overall evolution of the system is represented by an equilibrium of the corresponding MFG called Mean Field Equilibria (MFE). The discretization of the problems introduced by \cite{huang2006large} and \cite{lasry2007mean} has raised questions in numerical analysis and has inspired the study of MFG models in discrete time with a finite number of states, as seen in the work of \cite{gomes2010discrete}. These models are particularly useful for exploring the relationship between transportation research and recent advancements in MFG theory. The existing literature on MFGs in transportation has primarily focused on road traffic, congestion modeling, and dynamic routing, where the framework's scalability and equilibrium properties are well established. Recent studies have extended MFG formulations to incorporate congestion effects, multi-class agents, and time-dependent routing choices, demonstrating their potential for large-scale mobility systems. However, direct applications of MFG in maritime transportation remain limited.\\
\indent This paper aims to develop a discrete- time and state MFG model that analyzes competition among seaports. There are two main perspectives from which this competition can be examined. The first perspective is \emph{microscopic}, focusing on individual ships and their decision-making processes based on personal attributes, such as the type of cargo they carry, as well as network-related factors, like the travel time between their current port and destination port. The second perspective is \emph{macroscopic}, which looks at the overall behavior of ships. The model proposed in this paper seeks to balance these microscopic and macroscopic approaches. Analyzing the situation solely from the perspective of a single ship is unrealistic, as shipping companies typically own multiple vessels. For instance, there are 7,327 shipping companies listed, complete with detailed information about their suppliers, buyers, vessels, and seaports (see, e.g., \url{https://panjiva.com/sitemap/scac_profiles_directory}). One example, the shipping company \texttt{Aa Freight Inc}, operates 820 vessels. Therefore, in practice, one should consider a limited number of large companies. We will model one representative company at each seaport, reflecting the collective decisions made by all companies with ships at that port. The objective function takes into account \emph{transportation costs}, \emph{congestion costs}, and the so-called \emph{commercial margin} -- the difference between the value of a certain good at the arrival port and its value at the departure port. A statistical methodology will be provided to demonstrate how to infer the model's parameters from real-world data. Furthermore, we will present and discuss the results of an empirical application of this methodology using the \texttt{ShipFix} dataset. 
\subsection{Literature Review}\label{subsection::literaturereview}
We will first review the literature on MFGs and transportation research, then investigate a series of studies that link game theory with seaports competition. The research gap regarding the application of MFGs in seaport competition will be identified.\\
\indent MFE is a game-theoretic framework designed to model complex dynamics that arise from the interactions of a large population of rational agents who seek to optimize their utility. The behavior of these agents is characterized by optimal control problems, and the impact of each individual agent on the overall system is considered to be negligible. For a comprehensive presentation of MFG theory and its applications from both analytic and probabilistic perspectives, see the works of \cite{huang2006large}, \cite{lasry2007mean}, \cite{cardaliaguet2010notes}, and \cite{carmona2018probabilistic}. By taking advantage of the symmetry present in the game and allowing the number of agents to approach infinity, we can simplify the problem. This simplification leads to a scenario in which the original game structure is preserved but can be effectively reduced to a problem involving a ``representative player." This representative player responds optimally to the average behavior of the entire population. Once the limit problem for this representative agent is solved, its solution can be applied to the large-population game, yielding an approximate Nash equilibrium with diminishing error as the number of agents increases. The solution for the representative agent is defined by a pair of equations: the backward Hamilton-Jacobi-Bellman (HJB) equation and the forward Fokker-Planck-Kolmogorov (FPK) equation. The HJB equation ensures optimality for each individual player, while the FPK equation guarantees the time consistency of the solution. MFG is generally challenging to solve due its forward-backward structure, and different numerical methods for solving MFE have been proposed (see, e.g., \cite{achdou2020mean} for a recent review).\\
\indent Mean Field Games (MFG) have proven beneficial for modeling dynamic decision-making processes involving multiple agents and have gained popularity in various fields, such as finance (e.g., \cite{gueant2010mean}), the transition to green energy (e.g., \cite{flandoli2025structural}), crowd motion (e.g., \cite{lachapelle2011mean}), and autonomous driving (e.g., \cite{mo2024game}).\\
\indent In the context of the connection between discrete-time and state MFGs and transportation research, \cite{baillon2008markovian} introduced a Markov framework suitable for traffic assignment problems, akin to the formulation in \cite{tanaka2020linearly}. Meanwhile, \cite{chevalier2015micro} modeled interactions among drivers on a straight road as a non-cooperative game, characterizing its mean-field equilibrium. Additionally, \cite{gueant2015existence} provided results on the existence and uniqueness of mean-field equations on graphs, which account for ``congestion effects." \cite{bauso2016density} examined a continuous-time Markov chain to model the aggregated behavior of drivers within a traffic network.\\
\indent In the previous models, a large number of identical agents can exist across a finite number of states, with each agent behaving individually and rationally, transitioning from one state to another based on specific optimal criteria. Decisions are guided solely by information accessible to all agents, which includes the current state and the fraction of agents in each state. Similar to non-cooperative games, interactions may occur between players across different states. Moreover, due to the substantial number of agents, the mean-field hypothesis posits that the fraction of players in each state at any given time is the critical information for understanding global evolution.\\
\indent In the context of a game-theoretic framework applied to seaport competition, \cite{anderson2008game} analyzed capacity investment decisions among the seaports of Busan, Korea, and Shanghai using a \emph{Bertrand competition model}. \cite{kaselimi2011game} examined a two-stage game in which terminal operators compete on both prices and output, considering their capacity through a \emph{Cournot-type competition} approach. Meanwhile, \cite{luo2012post} investigated a two-stage game involving the seaports of Hong Kong and Shenzhen, where the seaports make capacity investment decisions in response to increasing market demand and differing service levels, while also considering pricing strategies, employing techniques from \emph{Bertrand competition} and \emph{Nash equilibrium}. The same objective and seaports were later examined in \cite{do2015application} via a two-person game model with uncertain demand and payoff \emph{Nash equilibrium}; here, seaports decide to invest under consideration that demand is uncertain or payoff is uncertain. \cite{ishii2013game} analysed a two-person game model with stochastic demand via \emph{Nash equilibrium} in which seaports make pricing decisions in the time of capacity investment. \cite{zhuang2014game} investigate service differentiation for seaports that handle containerized cargo and dry bulk cargo using a \emph{Stackelberg game} model. In this scenario, the leading port determines the output volumes for both types of cargo operations, while the follower port decides on its output volumes for container and bulk cargo operations. \cite{yip2014modeling} analyze a two-stage game in which seaports make decisions regarding terminal awards, and terminals set port charges while competing in terms of quantity. Their study focuses on both inter-port and intra-port competition related to terminal concession awarding. \cite{nguyen2016analysis} examine a two-stage game utilizing Nash equilibrium to explore pricing decisions made by seaports aimed at maximizing profits. This analysis identifies network links among seaports in three Australian regions: Queensland, South Australia, and Victoria, along with Western Australia, and considers the strategic interactions involved. Lastly, \cite{ignatius2018cooperate} investigate whether an alliance should be formed between three leading transshipment seaports located in Malaysia (Port Klang and the Port of Tanjung Pelepas) and Singapore (Port of Singapore) through \textit{Cournot competition}. For further insights into this area, we refer readers to the review by \cite{pujats2020review}.\\
\indent However, to the best of our knowledge, MFGs in seaport competition is still a nascent field that entails a lot of open questions and challenges. We believe this work represents a significant step toward a more refined application of MFGs in maritime traffic management.
\subsection{Contribution of this Paper}\label{subsection::contribution}
This paper aims to develop a game-theoretic framework for maritime flow modelling using mean-field approximation. We will start from identifying the factors affecting shipping companies decisions on the selection of seaports and we propose a MFG model at a mesoscopic scale to describe seaport competition. Building on the MFG formulation, we will propose a statistical methodology to infer the parameters of the corresponding MFE from real-world data. For real-world validation, we analyze the \texttt{ShipFix} dataset of daily coal exports worldwide from 2015 to 2025.\\
\indent In particular, the contribution of this paper include the following:
\begin{itemize}
    \item  Propose a MFG model to understand maritime flows and the sources of decision-making. This model will consider transportation costs, expected profit margins from loading specific goods at seaports and unloading them at various destinations, and a congestion factor that represents the cost of accessing the destination port.
    \item Provide a clear solution to the stationary version of the MFG model by making certain approximations. Specifically, we assume quadratic costs and that a sufficient number of ships are involved, so that the stochastic transition matrix between seaports closely approximates its deterministic equivalent. Interestingly, the solution indicates that for each category of goods, the flow of transitions is determined by the first eigenvector of the optimal controls. These eigenvectors satisfy a system of linear equations that share the same mean field. Additionally, we present a condition that combines all the cost parameters to ensure the existence of this unique solution.
    \item Develop a statistical methodology to infer the parameters of the game, i.e., the costs and commercial margins' components, from the \texttt{ShipFix} dataset. The methodology does not require the application of, e.g., fictitious play techniques to stabilize its solution. We analyze the results to clarify how the various cost components are influenced by the empirical flows.
\end{itemize}
The remainder of this paper is organized as follows: Section \ref{section::problem_formulation_and_solution} presents the problem that this paper addresses and provides a solution based on certain approximations. In Section \ref{section::statistical_inference}, we develop the statistical methodology needed to infer the parameters of the MFE from real-world data. Section \ref{section::empirical_example} focuses on validating the procedure outlined in Section \ref{section::statistical_inference} using the \texttt{ShipFix} dataset. Finally, Section \ref{section::conclusion_and_future_research} concludes the paper and discusses future research directions.
\section{Problem statement and solution}\label{section::problem_formulation_and_solution}
We consider a platoon of seaports indexed by \( i \in \{1, 2, \ldots, K\} \), where \( K \) represents the total number of seaports, and a set of \( N \) goods, indexed by \( n \in \{1, 2, \ldots, N\} \). Each port \( i \) is home to a number of ships. Each ship aims to be loaded at a specific time \( t \). The companies that own the ships are responsible for deciding what fraction of their fleet should be dispatched from port \( i \) to port \( j \). We first make the following assumption for modelling MFGs:
\begin{assumption}\label{assumption::idiosyncrasy}
    \item Companies do not have idiosyncrasy.
\end{assumption}
This implies that, for instance, the \texttt{Mediterranean Shipping Company} would make the same decisions as \texttt{CMA CGM} if both companies operated identical ships in the same locations. Therefore, we can combine the decisions of all the companies at each port \( i \) to create a $K$-dimensional vector representing the desired transitions for the entire fleet present at time \( t \) in seaport \( i \), directing them to each potential seaport destination \( j \). We denote by $\Phi_{i,j}^n(t^+)$ the flow of ships from seaport $i$ to seaport $j$, carrying good $n$; the notation $t^+$ emphasizes that this decision is made after the fleet's capacity at time $t$ (i.e., the number of ships) is known. Each company aims to make optimal decisions by maximizing its driving cost function, which is characterized by the following factors; more details will be given in Subsection \ref{subsec::maritme_flow_policy}. 
\begin{itemize}
    \item   \emph{Commercial margins}, are affected by taxes such as tariffs, and docking costs (e.g., \cite{wiegmans2008port}). The latter are often penalyzed by a congestion factor, which can be applied either explicitly or implicitly through an auction system. This practice is commonly referred to as competition-based pricing (e.g., \cite{hunt2025differences}). We denote by $M_{i,j}^{n}(t)$ the profit margin each company can achieve by transporting a unit of good $n$ from seaport $i$ to seaport $j$.
    \item  \emph{Transportation costs}, which are assumed to be proportional to the distance or time required for transport from seaport \(i\) to seaport \(j\). There are Python packages available, such as \texttt{Searoutes py} (\cite{searoute}), that can generate the shortest sea route between two points on Earth. We define transportation costs using the following relation: \(c_n g(T_{i,j})\), where \(c_n\) is a positive constant that depends on the type of good, and \(g:\mathbb{R}_{+}\rightarrow\mathbb{R}_{+}\) is a function, which is typically assumed to follow a power law based on the number of units to be transported (e.g., \cite{micco2002determinants}).
    \item \emph{Congestion costs}, which are affected by the flow arriving at destination ports \( j \). This aspect is more complex to manage than transportation costs because the flow \( \Phi_{i,j}^n(t^+) \) will reach port \( i \) at time \( t^{+} + T_{i,j} \). Therefore, the expected crowd at port \( j \) upon arrival can be expressed as:
\[
\sum_{\ell \neq i} \Phi_{\ell,j}^n(t^{+} + T_{i,j} - T_{\ell,j}).
\] 
\noindent This equation accounts for the flow from all other ports \( \ell \) that will arrive at port \( j \) before the flow from port \( i \).
\end{itemize}

In this paper, we aim to address the \textit{stationary} version of the MFG problem. Specifically, under the previous assumption, the flows \(\Phi_{i,j}^{n}(t^{+})\) exhibit the same distribution regardless of the time considered. This means that, at the appropriate time scale, quantities like \(\sum_{\ell \neq i} \Phi_{\ell,j}^n\) are treated as random variables rather than stochastic processes, which simplifies the formulation of congestion costs, for example. In particular, We need the following assumption.
\begin{assumption}\label{assumption::stationarity}
Let \( F^{n} \) be the total transportation capacity for a given good \( n \in \{1,\ldots,N\} \). We assume that \( F^{n} \) remains constant, meaning that companies are neither adding nor removing ships from the ecosystem we are considering.
\end{assumption}
\subsection{Maritime flow MFG}\label{subsec::maritme_flow_policy}
In this subsection, we outline the problem statement for the stationary version of the maritime flow MFG, which is the primary focus of this paper. Specifically, we define the payoff and cost functions for the coordinator at each of the seaports under consideration that determines the routes taken by ships transporting goods between seaports. At any given time \( t \), the coordinator at seaport \( i \), \( i \in \{1,\ldots,K\} \), aims to select, for every good $n \in \{1,\ldots,N\}$ the \( K \times K \) transition matrix \(Q_{i}= [ Q_{i,j}^{n}(t^{+})]_{j,n}\). The entry $(i,j)$ of this matrix represents the fraction of ships waiting to be loaded with good \( n \) at port \( i \) that will be sent to port \( j \). Denoting by $P^{n}=[P^{n}_{i}(t)]_{i}$, $i \in \{1,\ldots,K\}$, the $K$-dimensional probability distribution of the location of the good $n$ at time $t$, the expected flow of good $n$ from the seaport $i$ to the seaport $j$ is given by 
\begin{equation}\label{eq::realized_flow}
    \Phi_{i,j}^{n,Q}(t^+) = F^n P^n_i(t)\, Q^n_{i,j}(t^+)    
\end{equation}
It's essential to understand that the previous relation does not involve randomness; $Q^n_{i,j}(t^+)$ represents the desired probability of sending ships from port $i$ to port $j$ that are waiting to transport good $n$. In particular, in this paper, we model decisions for large ports with intense ship traffic, assuming that the expected flow is realized. In other words, the random variations in concentration, as explained by the central limit theorem, are negligible enough to disregard. While examining these variations falls outside the scope of this paper, it is worth noting that extending our results to a stochastic version of the model is feasible. However, this would necessitate the introduction of convexity terms associated with cost, which are non-linear functions of ship flows. In the case of ports experiencing high traffic volumes, these terms are sufficiently small that they are unlikely to significantly alter the simplicity and clarity of our results. Our findings maintain their robustness, making them both impactful and easy to understand.\\
\indent We now detail the cost functional. The coordinator in the seaport $i$ aims to select, for every good $n$, the optimal transition vector $Q_{i,\cdot}^{n}(t^{+})$ by minimizing her/his expected cost component of the cost functional, defined by the following two components:
\begin{itemize}
    \item Expected transportation cost at the decision time: 
    \begin{equation}\label{eq::expected_transportation_cost}
        C_{i,\cdot}^{n,Q} = c_n \sum_{\substack{j=1 \\ j\neq i}}^{K} \left|\Phi_{i,j}^{n,Q}(t^{+})\right|^{\gamma} g(T_{i,j}),
    \end{equation}
    where all the quantities have been previously introduced and $\gamma$ is a positive constant. Setting \(\gamma = 2\) simplifies the mathematical formulation. 
    \item The cost of congestion, or crowding, necessitates an understanding of its structure. To comprehend this cost, we must recognize that the coordinator decides to send goods \(n\) to seaport \(j\) while considering all types of goods arriving at \(j\) simultaneously. Specifically, the goods sent from seaport \(i\) at time \(t^{+}\) will arrive at seaport \(j\) at time \(\tau = t^{+} + T_{i,j}\). Similarly, goods sent from any other port \(\ell \neq j\) to the same seaport \(j\) and arriving at the same time \(\tau\) must be sent at time \(\tau - T_{\ell,j}\), which can be expressed as \(t^{+} + T_{i,j} - T_{\ell,j}\). By using the assumption that the expected flow is realized (see above), the cost of congestion can be written in the following way:
    \begin{equation}\label{eq::cost_of_congestion}
        R^Q_i(t^+)= \sum_{j=1}^{K} r_j \left| \sum_{n=1}^{N} \sum_{\substack{\ell=1 \\ \ell\neq i}}^{K} \Phi^{n,Q}_{\ell,j}(t^++T_{i,j}-T_{\ell,j}) + \sum_{n=1}^{N} \Phi^{n,Q}_{i,j}(t^+)\right|^\rho, 
    \end{equation}
    where $\rho$ is a positive constant. Again, setting \(\rho = 2\) simplifies the mathematical formulation. In particular, the cost of congestion described in Equation \eqref{eq::cost_of_congestion} is influenced not only by the control actions taken by the coordinator at seaport \(i\), but also by the decisions made at all other seaports \(\ell \neq i\), which are beyond the coordinator's control. The term $\sum_{n=1}^{N} \sum_{\substack{\ell=1 \\ \ell\neq i}}^{K} \Phi^{n,Q}_{\ell,j}(t^++T_{i,j}-T_{\ell,j})$ is therefore the so-called \emph{mean-field} term. Notice that the assumption that the expected flow is realized entails that the travelling times (or distances) are not random or that their randomness can be reduced to a change in the term $r_j$. As said before, in this paper, we aim to address the stationary version of the MFG problem, and therefore the cost of congestion reads as
        \begin{equation}\label{eq::cost_of_congestion_stationary}
        R^Q_i(t^+)= \sum_{j=1}^{K} r_j \left| \sum_{n=1}^{N} \sum_{\substack{\ell=1 \\ \ell\neq i}}^{K} \Phi^{n,Q}_{\ell,j} + \sum_{n=1}^{N} \Phi^{n,Q}_{i,j} \right|^\rho, 
    \end{equation}
\end{itemize}

\indent On the other hand, we define the expected coordinator's gain as the difference between the value of good $n$ at the destination seaport $j$, denoted as $v_{j}^{n}$, and its value at the departure seaport $i$, represented as $v_{i}^{n}$. We have
\begin{equation}\label{eq::gain}
    \mathbf{M}^{n,Q}_{i} = \sum_{j=1}^{K} \Phi_{i,j}^{n,Q}(t^+)[v_j^n-v_i^n]:=\sum_{j=1}^{K} \Phi_{i,j}^{n,Q} M^n_{i,j},
\end{equation}
where in the last definition, we use the stationarity assumption and the notation $M^n_{i,j}$ to denote the commercial margin $v_j^n-v_i^n$.\\
\noindent Therefore, the cost functional that the coordinator at the seaport $i$ managing the good $n$ aims to maximize is defined as:
\begin{equation}\label{eq::cost_functional}
    J_{i,n}^{N}(Q_{i,\cdot}^{n}, Q_{-i,\cdot}) = \underbrace{\sum_{j=1}^{K} \Phi^{n,Q}_{i,j} M^n_{i,j}}_{\text{expected gain}} - \underbrace{\sum_{j=1}^{K} r_j \left| \sum_{\substack{m=1 \\ m \neq n}}^{N} \sum_{\substack{\ell=1}}^{K} \Phi^{m,Q}_{\ell,j} + \sum_{\substack{\ell=1 \\ \ell \neq i}}^{K} \Phi^{n,Q}_{\ell,j} + \Phi^{n,Q}_{i,j}  \right|^\rho}_{\text{cost of congestion}} - \underbrace{c_n \sum_{\substack{j=1 \\ j\neq i}}^{K} \left|\Phi_{i,j}^{n,Q} \right|^{\gamma} g(T_{i,j})}_{\text{transportation cost}}.
\end{equation}
This quantity depends not only on the control exerted by the coordinator at seaport \(i\), denoted as \(Q_{i,\cdot}\), but also on the policies of the other coordinators, \(Q_{-i,\cdot}\), at the other seaports. This relationship is represented through the following terms: 
\[
\sum_{\substack{m=1 \\ m \neq n}}^{N} \sum_{\ell=1}^{K} \Phi^{m,Q}_{\ell,j} + \sum_{\substack{\ell=1 \\ \ell \neq i}}^{K} \Phi^{n,Q}_{\ell,j}.
\] 
These terms illustrate how decisions made by coordinators at various seaports affect the overall outcome. We are dealing with a non-cooperative game. Specifically, under Assumption \ref{assumption::stationarity}, the mean-field is represented as a \(N \times K\) matrix \([\varphi]_{n,i} := \varphi_{i}^{n}\), where the \((n,i)\)-\emph{th} component indicates the quantity of good \(n\) available at seaport \(i\). In this context, the controlled flow \(\Phi^{n,Q}_{i,j}\) is defined as the product of the mean-field (the quantity at the seaport \(i\)) and the transition (the amount sent from this seaport to other seaports):
\begin{equation}\label{eq::mean_field}
    \Phi^{n,{Q}}_{i,j} = \varphi_{i}^{n} \cdot Q_{i,j}^{n}.
\end{equation}
\subsubsection{Solution of the quadratic version of the maritime flow MFG}
We consider the cost functional in Equation \eqref{eq::cost_functional} when $\rho=\gamma=2$ and we take into account Equation \eqref{eq::mean_field}. Whence, the optimization problem becomes
\begin{small}
\begin{equation}\label{eq::optimization_quadratic}
        \max_{\substack{Q_{i,\cdot}^{n}}}\left\{\sum_{j=1}^{K} \varphi_{i}^{n} \cdot Q_{i,j}^{n} M^n_{i,j} - \sum_{j=1}^{K} r_j \left[ \sum_{\substack{m=1 \\ m \neq n}}^{N} \sum_{\substack{\ell=1}}^{K} \varphi_{\ell}^{m} \cdot Q_{\ell,j}^{m} + \sum_{\substack{\ell=1 \\ \ell \neq i}}^{K} \varphi_{\ell}^{n} \cdot Q_{\ell,j}^{n} + \varphi_{i}^{n} \cdot Q_{i,j}^{n}
        \right]^2 - c_n \sum_{\substack{j=1 \\ j\neq i}}^{K} (\varphi_{i}^{n} \cdot Q_{i,j}^{n})^2 g(T_{i,j})\right\},
        \end{equation}
\end{small}
coupled with the following static mass conservation equation:
\begin{equation}\label{eq:FP:MF}
    \forall n \in \{1,\ldots,N\}:\:\left\{ \forall j \in \{1,\ldots,K\}\,:\, \vphjnd = \sum_{i=1}^{K} \vphind \Qcont\right\}.
\end{equation}
The previous equation is saying $\textcolor{black}{\varphi^n}$ is the stationary distribution for $\textcolor{black}{Q^n}$.\\
\indent We are now ready to state the proposition that derives the mean-field approximation on ports of origin.
\begin{proposition}\label{prop::mean_field_approximation}
consider the constrained maximization problem in Equations \eqref{eq::optimization_quadratic}-\eqref{eq:FP:MF}. If the number of seaports \(N\) is sufficiently large, each coordinator focuses on optimizing the problem regarding one seaport of origin and one good, treating the other flows as independent effects:
\begin{equation}\label{sim:mfg:limit}
\sum_{\substack{\ell=1 \\ \ell \neq i}}^{K} \varphi_{\ell}^{n} \cdot Q_{\ell,j}^{n} + \varphi_{i}^{n} \cdot Q_{i,j}^{n}  \simeq \sum_{\ell=1}^{K} \varphi^n_\ell Q^{n}_{\ell,j} + \vphind \Qcont = \varphi^n_j + \vphind \Qcont.
\end{equation}
Therefore, the maximization problem in Equation \eqref{eq::optimization_quadratic} reads as:
\begin{equation}\label{eq:max:quad:MF:red}
    \max_{\substack{Q_{i,\cdot}^{n}}}\left\{\sum_{j=1}^{K} \varphi_{i}^{n} \cdot Q_{i,j}^{n} M^n_{i,j} - \sum_{j=1}^{K} r_j \left[ \sum_{\substack{m=1 \\ m \neq n}}^{N}  \varphi_{j}^{m}  +  \varphi_{j}^{n} \cdot Q_{\ell,j}^{n} + \varphi_{i}^{n} \cdot Q_{i,j}^{n}
        \right]^2 - c_n \sum_{\substack{j=1 \\ j\neq i}}^{K} (\varphi_{i}^{n} \cdot Q_{i,j}^{n})^2 g(T_{i,j})\right\},
\end{equation}
\end{proposition}
\proof{Proof of Proposition \ref{prop::mean_field_approximation}} 
The approximation in Equation \eqref{sim:mfg:limit} is the typical simplification done in the MFG framework, namely one replaces multiple interactions with a single interaction involving the mean-field and one ``isolated control." In this specific case, it is indicated that when a coordinator focuses on controlling \( Q^n_{i,j} \), they are not directly aware of the isolated decisions made at other ports of origin. From a continuum perspective of ports, this can be expressed as:

\[
\int_{\ell\neq i} \varphi^n_\ell Q^{n}_{\ell,j} \, dm(\ell) = \varphi^n_j.
\]    
Under the previous approximation, the  the equation \( \sum_{\ell} \varphi^n_\ell Q^{n}_{\ell,j} = \varphi^n_j \) simplifies to Equation \eqref{eq:FP:MF}), which represents the stationarity of the mean-field.
\endproof
\vspace{0.5cm}
Notice that Equations \eqref{eq:max:quad:MF:red} are coupled via the $N \times K$-dimensional mean-field, updated via Equation \eqref{eq:FP:MF}. We now provide the solution of the quadratic version of the maritime flow MFG.
\begin{theorem}\label{theorem::MFG_solution_quadratic}
   Consider the constrained maximization problem in Equations \eqref{eq:max:quad:MF:red}-\eqref{eq:FP:MF}. For every origin seaport $i \in \{1,\ldots,N\}$, every destination port $j \in \{1,\ldots,N\}\setminus\{i\}$  and for every good $n \in \{1,\ldots,N\}$, let 
   \begin{equation}\label{eq:quad:names}
  \wM_{i,j}^n:=\frac{M^n_{i,j}}{2 },\quad w^n_{i,j}=\frac{1}{r_j+c_n\, g(T_{i,j})\one_{\{j\neq i\}}}, \quad 
 \uw_{i,j}^n:=\frac{w_{i,j}^n}{\sum_{\substack{\ell=1 \\ }}^{K} w_{i,\ell}^n},\quad \vphjdd=\sum_{m=1}^{N}\vphjmd
\end{equation} 
where $\one_{\{j\neq i\}}$ indicates the case in which the seaport $i$ is different from the seaport $j$.   Then the solution of \eqref{eq:max:quad:MF:red}-\eqref{eq:FP:MF} is given by 
\begin{equation}\label{eq:quad:Q}
    Q^n_{i,j} = \frac{w^n_{i,j}}{\varphi^n_i} \left\{ \displaystyle \wM_{i,j}^n - \sum_{\substack{\ell=1 \\ }}^{K} \uw_{i,\ell}^n \,\wM_{i,\ell}^n - \left[ r_j \vfb_j - \sum_{\substack{\ell=1 \\ }}^{K} \uw_{i,\ell}^n \,r_\ell \vfb_\ell\right]\right\} + \uw^n_{i,j}.
\end{equation}
\end{theorem}
\proof{Proof of Theorem \ref{theorem::MFG_solution_quadratic}} 
Let $i \in \{1,\ldots,K\}$ and $n \in \{1,\ldots,N\}$ be a fixed seaport and a fixed good. Let $[q_j]_{j}$ be the $K$-row of the transition matrix $Q_i$. This vector reflects the fraction of the available good at port \(i\) for good \(n\) that will be transferred to port \(j\). It's worth noting that the configuration where \(j = i\) is allowed, representing the portion of the good that remains at port \(i\). Now, let $\lambda_{i}^{n} \in \mathbb{R}$ be a Lagrange multiplier to keep $[q_j]_{j}$ in the simplex. Then, the (augmented) first order condition is the following
\begin{equation*}
-\lambda \sum_{j=1}^{K} q_j + \vphind \sum_{j=1}^{K} q_j M^n_{i,j} - \sum_{j=1}^{K} r_j \left[ \left(\sum_{\substack{m=1 \\ m\neq n}}^{n} \vphjmd \right) + \vphjnd + \vphind q_j  \right]^2 - c_n \sum_{\substack{j=1 \\ j \neq i}}^{K} \left(\vphind q_j\right)^2\cdot g(T_{i,j})=0,
\end{equation*}
from which we obtain, for every $i \in \{1,\ldots,K\}$ and $n \in \{1,\ldots,N\}$, that 
\begin{equation*}
    \frac{\lambda}{\varphi^n_i}=:\lambda_i^n = M^n_{i,j} - 2 r_j (\vphjdd + \vphind q_j) - 2 c_n g(T_{i,j}) \vphind q_j \one_{\{j\neq i\}}.
\end{equation*}
At this point, we write each $[q_j]_{j}$ -- i.e., $Q_{i,j}^{n}$ for fixed $i$ and $n$ -- as a function of the Lagrange multiplier  
    \begin{equation}\label{eq:quad:sol}
    Q_{i,j}^{n} = \frac{w^n_{i,j}}{2\vphind} \left(M^n_{i,j}-2r_j \vphjdd -\lambda_i^n\right),
    \end{equation}
where the Lagrange multiplier is obtained via saturating the constraint $\sum_{j=1}^{K} \Qcont = 1$. We obtain that 
\begin{equation*}
     \lambda_i^n = \frac{1}{\sum_{j=1}^{K} w_{i,\ell}^n}\left( \sum_{\ell=1}^{K} w_{i,\ell}^n(M_{i,\ell}^n -2 r_{\ell}  \vphldd )\right) - 2 \frac{1}{\sum_{\ell=1}^{K} w_{i,\ell}^n}\varphi_i^n.
\end{equation*}
By substituting the previous expression into Equation \eqref{eq:quad:sol}, we obtain Equation \eqref{eq:quad:Q}. This ends the proof.
\endproof
\vspace{0.5cm}
\indent In particular, in order to obtain the mean-field term, one can compute the first eingenvector of the transition matrix $[Q_{i,j}^{n}]_{i,j}$ for good $n$, that is $\varphi_{i}^{n}$. In particular, this can be done for every good $n \in \{1,\ldots,N\}$. Whence, we can revisit the optimal control problem until we reach a fixed point.

\indent Before proceeding, we now make the following remarks.
\begin{remark}
We would like to pinpoint the interpretability, for every $i \in \{1,\ldots,K\}$ and $n \in \{1,\ldots,N\}$, of Equation \eqref{eq:quad:Q}. First, as the expected margins $\wM_{i,j}^n$ increase -- relative to the barycenter of the margins at destination ports, represented by $\sum_{\ell=1}^{K} \uw_{i,\ell}^n \,\wM_{i,\ell}^n$ -- the quantity dispatched to port $j$ significantly rises. Second, the greater the anticipated crowding at the destination -- again, relative to the barycenter of the crowd at destination ports -- the smaller the quantity that will be dispatched there. Third, the overall costs, particularly transportation costs, are negatively affecting the quantity being sent. Finally, the averaging scheme assigns less weight to locations where high costs are incurred (regardless of the items being sent, due to the term \( r_j \)). It also disregards places experiencing high transportation costs (because of the term \( c_n g(T_{i,j}) \)), with the exception of the designated port \( i \), where transportation costs are considered to be zero. Notice that the rems $Q_{i,j}^{n}$ and $\uw_{i,j}^{n}$ are dimensionless, $w_{i,j}^{n}$ is measured in cost's inverse, and $M_{i,j}^{n}$ must have the same units as $\varphi_{i}^{n}$ for the formula in Equation \eqref{eq:quad:Q}) to be valid.

\end{remark}

\begin{remark}
We do not address the positivity constraint in the optimization problem, as it can be quite challenging due to its potential to create sparsity in the solutions. However, it is clear that in realistic scenarios, there should exist a non-empty domain for the parameters of the problem (such as costs and commercial margins) that would prevent transitions from becoming negative. The reader may be convinced that negative transitions should not occur in practical situations. If it is more optimal to import goods from port \(j\) rather than export goods to port \(j\), then, within this framework, there is no reason for the export transition matrix to contain a negative term. It is sufficient for the export matrix, which originates from port \(j\), to have a positive term leading to port \(i\).
\end{remark}
\subsubsection{Derivation of the MFE for the maritime flow MFG}
If the number of goods approaches infinity \( N \rightarrow \infty \), we recover the problem for a representative good. In this case, one can express the mean-field \([\varphi_{i}]_{i}\) as the row eigenvector of \([Q_{,j}]_{i,j}\) associated with the eigenvalue of one. Importantly, the expression between the curly brackets in Equation \eqref{eq:max:quad:MF:red} is quadratic in $Q_{i,j}$, and it is such that the term $Q_{i,j}$ is always multiplied by $\varphi_{i}$; in this way, we have that the right-hand side of Equation \eqref{eq:quad:Q} has a prefactor $\varphi_{i}$ which cancels out when writing $\sum_{i=1}^{K}\varphi_{i}Q_{i,j}$. Precisely, we have that the following computations hold true:
\begin{eqnarray*}
\sum_{i=1}^{K} \varphi_i Q_{i,j}&=&\sum_{i=1}^{K} w_{i,j}\left[ \wM_{i,j} - \sum_{\ell=1}^{K} \uw_{i,\ell} \wM_{i,\ell} - \left( r_j \varphi_{j} - \sum_{\ell=1}^{K} \uw_{i,\ell} \,r_\ell\,\varphi_{\ell} \right)\right]\ + \sum_{i=1}^{K} \uw_{i,j}\varphi_i \\
\Rightarrow \left[ 1+ r_j \sum_{i=1}^{K} w_{i,j}\right]\cdot\varphi_j &=& \underbrace{\sum_{i=1}^{K} w_{i,j} \left[\wM_{i,j} - \sum_{\ell=1}^{K} \uw_{i,\ell} \wM_{i,\ell} \right]}_{\left( \sum_{i=1}^{K} w_{i,j}\right)\, \bM_j} + \sum_{\ell=1}^{K} \left(\sum_{i=1}^{K} w_{i,j} \uw_{i,\ell}\right) r_\ell \,\varphi_\ell + \sum_{i=1}^{K} \varphi_i \uw_{i,j}\\
\Rightarrow \left[ 1 + r_j \sum_{i=1}^{K} w_{i,j}\right]\cdot\varphi_j &=& \left(\sum_{i=1}^{K} w_{i,j}\right)\, \bM_j + \sum_{\ell=1}^{K} \left(\sum_{i=1}^{K} w_{i,j} \uw_{i,\ell}\right) r_\ell \,\varphi_\ell + \sum_{\ell=1}^{K} \varphi_\ell\, \uw_{\ell,j},
\end{eqnarray*}
where we use the notation  $\bM_j$ for the \emph{overall average relative margin}:
$$\bM_j=\sum_{i=1}^{K} \frac{w_{i,j}}{\sum_{i'=1}^{K} w_{i',j}} \left[\wM_{i,j} - \sum_{\ell=1}^{K} \uw_{i,\ell} \wM_{i,\ell} \right].$$
\indent We can write now the following theorem, which ensures the existence of a unique solution for the MFE of the maritime flow MFG.
\begin{theorem}\label{thm::existence_of_the_equilibrium}
Consider the MFE problem defined by the following two equations
\begin{small}
\begin{equation*}
    Q_{i,j} = \frac{w^n_{i,j}}{\varphi_i} \left\{ \displaystyle \wM_{i,j} - \sum_{\substack{\ell=1 \\ }}^{K} \uw_{i,\ell} \,\wM_{i,\ell} - \left[ r_j \varphi_j - \sum_{\substack{\ell=1 \\ }}^{K} \uw_{i,\ell} \,r_\ell \varphi_\ell\right]\right\} + \uw^n_{i,j},\quad \forall j \in \{1,\ldots,K\}\,:\,\varphi_{j} = \sum_{i=1}^{K}\varphi_i Q_{i,j},
\end{equation*}
where the first equation defines the control $Q_{i,j}$, i.e., the probability of sending the representative good from the seaport $i$ to the seaport $i$, and the second equation characterizes the stationarity of the mean-field. Then, the MFE problem has a unique solution if and only if  
   \begin{equation}\label{eq:lin:det}
    \det\left( \left[ \left( 1 + r_j \sum_{\substack{\ell=1 \\ }}^{K}  w_{i,j} \right) \one_{\{j=\ell\}} - \left\{ \left(\sum_{\substack{\ell=1 \\ }}^{K}  w_{i,j} \uw_{i,\ell}\right) r_\ell  + \uw_{\ell,j} \right\}
    \right]_{j,\ell} \right) \neq 0, 
    \end{equation}  
\end{small}
where, in general, $\det(A)$ denotes the determinant of the matrix $A$.
\end{theorem}
\proof{Proof of Theorem \ref{thm::existence_of_the_equilibrium}} 
The proof is a continuation of the above computations. In particular, by using the notations:
\[
    R_j := {1\over r_j}+ \sum_{\substack{i=1 \\ }}^{K} w_{i,j},\quad %
    m_j := \frac{\sum_{\substack{i=1 \\ }}^{K} w_{i,j}}{r_j}\, \bM_j,\quad %
    C_{j,\ell} := \left(\sum_{\substack{i=1 \\ }}^{K} w_{i,j} \uw_{i,\ell}\right) \frac{r_\ell}{r_j}  + \frac{\uw_{\ell,j}}{r_j},
\]
the first equation of the MFE problem reads as $R_j \varphi_j = m_j + \sum_{\ell=1}^{K} C_{j,\ell} \, \varphi_\ell$, which gives a \emph{linear expression} for $\varphi$:
\begin{equation}\label{eq:phi:lin:sum}
        \sum_{\ell=1}^{K} \left[ R_j \one_{\{j=\ell\}} - C_{j,\ell} \right] \cdot \varphi_\ell = m_j.
\end{equation}
The previous equation can be rewritten, in algebraic notations, as 
    \begin{equation}\label{eq:phi:lin}
        \Omega \cdot \phi = m.
    \end{equation}
where $m$ is the vector of the $m_j$ components and $\Omega$ a matrix with the components $\Omega_{j,\ell}:=R_j \one_{\{j=\ell\}} - C_{\ell,j}$. In other terms, the fixed point exists if and only if $\Omega$ is an invertible matrix. In particular, a way to write this condition is that the upper determinant is not null. This ends the proof.
\endproof
\vspace{0.5 cm}
\indent Notice that the existence of a MFE is solely contingent on the cost structure, while the control variable \(Q_{i,j}\) also depends on profit margins. More specifically, the determinant condition outlined in Equation \eqref{eq:lin:det} relies exclusively on the cost parameters \(\uw_{i,j}\) -- which are influenced by \(r_{j}\) and \(g(T_{i,j})\) -- as well as the cost  \(c_n\). These parameters represent the direct transportation costs, interaction terms that reflect cross-route frictions or complementarities, and the return associated with each destination, respectively. Therefore, the existence and uniqueness of the MFE are determined entirely by the structure of these costs. This finding suggests that, economically, the cost structure establishes the geometry of feasible trade flows across the network of seaports. It governs whether the decentralized decisions of individual coordinators can combine to form a globally consistent stationary distribution. When the cost matrix is well-conditioned -- i.e., when the determinant in Equation \eqref{eq:lin:det} is non-zero -- the mapping from local shipping probabilities \(Q_{i,j}\) to the aggregate distribution \(\varphi\) becomes invertible. This ensures that the optimal routing decisions made by coordinators are mutually consistent and lead to a unique stationary MFE. Conversely, if the cost structure makes this mapping singular, the system loses its ability to determine a unique equilibrium distribution. This situation can arise, for example, when transportation costs are perfectly symmetric or when specific combinations of routes create effective zero-cost cycles, rendering it indifferent for coordinators to ship through one seaport or another. In such cases, multiple MFEs may align with cooperators' incentives, or none may exist at all. The determinant condition, therefore, offers a precise analytical framework for determining when the equilibrium operator is well-defined, linking the mathematical requirement of invertibility to an economically intuitive concept of network regularity and the absence of degeneracies in the underlying cost structure. While the existence and uniqueness of the MFE rest solely on costs, the actual equilibrium pattern of flows \(Q_{i,j}\) is also influenced by the margins \(\wM_{i,j}\). These margins represent destination-specific profitability differentials or other factors that affect coordinators' routing incentives without altering the fundamental feasibility of the equilibrium. Once the equilibrium mapping is clearly defined, variations in margins shift the relative attractiveness of destinations, which in turn influences the composition of trade flows across the network.\\
\noindent In summary, while costs determine whether an equilibrium can exist at all, margins dictate where within the feasible set of flows the economy ultimately settles. In this context, costs define the structural stability of the equilibrium, while margins shape its economic content -- specifically, the distribution of activities among the seaports once consistency is established.
\section{Statistical inference of the MFG parameters}\label{section::statistical_inference}
In this section, we develop a statistical methodology to infer the parameters of the quadratic version of the maritime flow MFG from real-world data, whose solution is provided in Theorem \ref{theorem::MFG_solution_quadratic}.  We will validate this methodology in Section \ref{section::empirical_example}. Consider, again, Equation \eqref{eq:quad:Q}, opportunely rewritten:
\begin{equation}\label{eq:quad:Q:2}
    \varphi^n_i Q^n_{i,j} = w^n_{i,j}\left\{ \displaystyle \wM_{i,j}^n - \sum_{\substack{\ell=1 \\ }}^{K} \uw_{i,\ell}^n \,\wM_{i,\ell}^n - \left[ r_j \sum_{m=1}^{N}\vphjmd - \sum_{\substack{\ell=1 \\ }}^{K} \uw_{i,\ell}^n \,r_\ell \sum_{m=1}^{N}\vphjmd \right]\right\} + \uw^n_{i,j} \varphi^n_i.
\end{equation}
First, we notice that all the terms between curly brackets can be expressed as functions of the MFG's parameters; see Equation \eqref{eq:quad:names}. Precisely, the expected margins $\wM^n_{i,j}$ are defined as $2 \wM^n_{i,j}=M_{i,j}^{n}=v_{j}^{n}-v_{i}^{n}$, the weights $\uw_{i,j}^n$ are normalized version of $w_{i,j}^{n}$, defined as:
\begin{equation}
  \label{eq:weight:approx}
  w^n_{i,j}=\frac{1}{r_j+c_n\, g(T_{i,j})\one_{\{j\neq i\}}},
\end{equation}
which are direct function of the congestion or crowdedness coefficient $r_j$, the transportation cost coefficient $c_n$, and the distance $T_{i,j}$ between two seaports $i$ and $j$. Our approach will consist of two steps.
\begin{itemize}
    \item \textit{Step 1}: Reformulate Equation \eqref{eq:quad:Q:2} so that it can be addressed as a solution to a linear regression problem, utilizing the data available in our real-world database, which will be described in Section \ref{section::empirical_example}.
    \item \textit{Step 2}: Express \(v_{i}^{n}\), \(c_n\), and \(r_j\) in terms of the coefficients identified in the linear regression from the previous step. This formulation should enable their retrieval through a cost-effective minimization algorithm, such as those available in common Python libraries like \texttt{Scipy}.
\end{itemize}

\indent In order to deal with the first step, it is essential to notice that a real-world dataset can violate the assumption of stationarity. To account for the just-mentioned violation, we will introduce the following time dependency in Equation \eqref{eq:quad:Q:2}
\begin{eqnarray*}\label{eq:quad:Q:2:1}
    \locf{\varphi^n_i Q^n_{i,j}}_t = w^n_{i,j} \locf{\displaystyle \wM^n_{i,j} - \sum_{\ell=1}^{K} \uw^n_{i,\ell} \,\wM^n_{i,\ell}}_\infty + w^n_{i,j} \sum_{\ell=1}^{K} (\uw^n_{i,\ell} -\one_{\{\ell=j\}}) \,r^n_\ell  \locf{ \sum_{m=1}^{N} \varphi^m_\ell}_{t:}+ \uw^n_{i,j} \locf{\varphi^n_i }_{t},
\end{eqnarray*}
where the low brackets highlight the time dependency and the subindex indicates what is the exact time scale of the real-word dataset that we will consider. In detail we have that:
\begin{itemize}
    \item The quantity $\locf{\varphi^n_i Q^n_{i,j}}_t$ is considered at a daily frequency, and it represents the exports.
    \item  We consider the relative margin $\locf{\displaystyle \wM^n_{i,j} - \sum_{\ell=1}^{K} \uw^n_{i,\ell} \,\wM^n_{i,\ell}}_{\infty}$ as constant over the real-word dataset, and therefore we use the subindex $\infty$. Admittedly, this is a limitation of our statistical inference framework; we discuss future research direction in Section \ref{section::conclusion_and_future_research}.
    \item The mean field, specifically the occupation at destinations \( \locf{\varphi^m_j}_{t:} \) and \( \locf{\varphi^m_\ell}_{t:} \), which is included in the relative mean field component, should be interpreted as \emph{the expected occupation at ports \( j \) and \( \ell \) from the perspective of date \( t \)}. We denote it with a subindex \( t: \) because it would corresponds to the expected cost of congestion by considering a random transportation time. In our estimation framework, we will calculate the expected usage of a port by averaging the imports over several days, along with considering typical transportation times; again, this is an approximation that can be refined in future work.
    \item  The last term is clear: $\locf{\varphi^n_i }_{t}$ indicates the amount of good available at time $t$ in port $i$ to be shipped elsewhere. We consider imports at date $t$.
\end{itemize}

\indent With the previous time-dependency conventions defined, we now refine the two steps above.
\begin{itemize}
    \item \emph{Step 1: A linear regression.}   Let $A_{i,j}, B_{i,j,\ell}, C_{i,j}$ be the following quantities  
    \begin{equation*}
     A_{i,j} \simeq  w^n_{i,j} \locf{\displaystyle \wM^n_{i,j} - \sum_{\ell=1}^{K} \uw^n_{i,\ell} \,\wM^n_{i,\ell}}_\infty,\quad B_{i,j,\ell} \simeq  w^n_{i,j} \sum_{\ell=1}^{K} (\uw^n_{i,\ell} -\one_{\{\ell=j\}}) \,r^n_\ell,\quad C_{i,j}      \simeq \uw^n_{i,j},
    \end{equation*}
where the symbol $``\simeq"$ is symptomatic of the fact that we potentially do not account for the total time variability of the real-world data set, which is reflected by the error term in our resulting linear regression model here below:
\begin{equation}
  \label{eq:estimate:model}
  \locf{\varphi^n_i Q^n_{i,j}}_t = A_{i,j} + \sum_{\ell=1}^{K} B_{i,j,\ell} \locf{\varphi^m_\ell}_{t:} + C_{i} \locf{\varphi^n_i }_{t} + \varepsilon_{i,j}(t).
\end{equation}
\item \emph{Step 2: A simple optimization.} Once the linear regression(s) in Equation \eqref{eq:estimate:model} are performed, we are ready to estimate the MFG's parameters when observational data become available. To maintain linearity with respect $w_{i,j}^{n}$, we perform the following approximation 
\begin{equation*}
    A_{i,j}  \simeq w_{i,j}^{n} \locf{\displaystyle \wM^n_{i,j} - \frac{1}{K-1} \sum_{\ell=1}^{K} \wM^n_{i,\ell}}_\infty.
\end{equation*}
By substituting $\wM^n_{i,\ell} = \tilde{v}_j-\tilde{v}_i$, where $v_{k}=\tilde{v}_k-\frac{1}{K-1}\sum_{\ell=1}^{K}\tilde{v}_{\ell}$, and the definition of $w_{i,j}^{n}$, we obtain the following (for the sake of notation, we omit the symbol $\locf{\,\,\cdot\,\,}_{\infty}$):
\begin{equation*}
    A_{i,j}  \simeq w_{i,j}^{n} ( v_j-v_i ) = \frac{v_j-v_i}{r_j + c_n g(T_{i,j})\one_{\{\ell=j\}}}.
\end{equation*}
We solve the previous expression by minimizing over $c_n, (r_j, v_{j})_{1 \leq j \leq K}$, using a numerical solver, the following expression:
\begin{equation*}
   \sum_{\substack{i,j=1 \\ i \neq j}}^{K} \left\|A_{i,j} \cdot (r_j+c_ng(T_{i,j})\one_{\{\ell=j\}}) - (v_j-v_i) \right\|^{2}.
\end{equation*}
\end{itemize}
We summarise the just described two-step procedure in the following proposition.
\begin{proposition}\label{prop::two_steps_procedure}
    Consider the MFE problem's control $Q_{i,j}$ in Equation \eqref{eq:quad:Q:2} and the time-convention introduced in the previous paragraphs. Then the MFG's parameters $c_{n}$, $n \in \{1,\ldots,N\}$, $r_j$, $v_j$, $1 \leq j \leq K$, can be obtained using the following two-step procedure.
    \begin{itemize}
        \item \textbf{\emph{Step 1}}: For each pair of source and destination seaports $(i,j)$, obtain from observational data the incoming flows and transitions and define the following quantity:
        \begin{equation*}
            (Y,X_1,\ldots,X_K,Z)=\left( \locf{\varphi^n_i Q^n_{i,j}}_t, \locf{\varphi^m_1}_{t:},\ldots, \locf{\varphi^m_K}_{t:}, \locf{\varphi^n_i }_{t}\right)_{t}
        \end{equation*}
        Then, perform the $K(K-1)$ linear regressions corresponding to Equation \eqref{eq:estimate:model}.
        \item \textbf{\emph{Step 2}}: Insert the calculated $A_{i,j}$ values and the distances $T_{i,j}$ between sources and destinations in the following minimization problem
        \begin{equation*}\label{eq::minimization}
   \min_{\substack{c_n, (r_j,v_j)_{1 \leq j \leq K}}}\sum_{\substack{i,j=1 \\ i \neq j}}^{K} \left\|A_{i,j} \cdot (r_j+c_ng(T_{i,j})\one_{\{\ell=j\}}) - (v_j-v_i) \right\|^{2}.
\end{equation*}
    \end{itemize}
\end{proposition}
It is worth noticing that the previous stepwise minimization process can be improved in several ways. By performing both minimizations simultaneously, we can take advantage of the inherent nature of linear regression, which is all about minimizing the mean squared error (MSE). Additionally, using $\uw_{i,j}^n$ instead of $w_{i,j}^n$ is a worthwhile consideration, even though it may complicate Step 2 for the minimizer. Nevertheless, at this stage of the theory of MFG for maritime traffic, we are committed to maintaining a straightforward and comprehensible approach. We will pursue these enhancements in future work. 
\section{Real-world data application}\label{section::empirical_example}
We use the \texttt{ShipFix} dataset, where Shipfix is a company that utilizes a decentralized email exchange among ship owners through the \texttt{Veson proprietary SaaS} platform. This platform is designed for managing route reservations and other administrative processes related to global shipping operations. Over the past 10 years, billions of data points have been collected, anonymized, and aggregated across $243$ cargo types. The cargo type ``Dry Coal" represents the largest segment, accounting for approximately $15\%$ of the total data points, while the second-largest category only makes up $5.7\%$; see Table \ref{tab:shipfix:cargo}. 
\begin{table} 
  \centering\small\rowcolors{2}{tableShade}{white}
 \begin{tabular}{lr}
 & pct. of observations (\%) \\
Cargo type &  \\
\tt dry.coal & 15.36 \\
\tt dry.agri\_products.grains & 5.70 \\
\tt dry.agri\_products.grains.wheat & 4.63 \\
\tt dry.fertilizers & 4.15 \\
\tt dry.steel\_products & 4.01 \\
\tt dry.cement.clinker & 3.27 \\
\tt dry & 3.12 \\
\tt dry.steel\_products.steel\_coil & 2.63 \\
\tt dry.minerals & 2.55 \\
\tt dry.agri\_products.grains.corn & 2.54 \\
\tt dry.petcoke & 2.47 \\
\tt dry.fertilizers.urea & 2.39 \\
\end{tabular}
  \caption{The table shows the twelve most frequent cargo types in the main \texttt{ShipFix} dataset, expressed in percentages.}
  \label{tab:shipfix:cargo}
\end{table}

\indent We evaluate our statistical methodology in Section \ref{section::statistical_inference} concerning the cargo type ``Dry Coal," which is indeed the most commonly exchanged raw material between seaports. Additionally, we chose to conduct this statistical study at the country level rather than focusing on individual seaports. While this approach can be refined in future research, the methodology will remain unchanged. The results we will obtain at the country level are compelling enough to warrant more precise estimations in subsequent studies.
\subsection{Data processing}\label{subsec::data_processing}
Data are extracted from the main importing and exporting countries. Table\ref{tab:yearly:exports:imports} presents the average yearly exports and imports for the period from January 2, 2015, to May 17, 2025 for these countries.

\begin{table}
  \centering\small
  \resizebox{.45\columnwidth}{!}{%
    \begin{tabular}[t]{lrrrrr|}
      {\bf Imports} & India & China & Vietnam & Japan & S.K. \\
      Indonesia & 165,802 & 164,098 & 26,333 & 11,316 & 15,152 \\
      Australia & 82,115 & 31,210 & 19,276 & 19,470 & 8,412 \\
      SAFR & 52,190 & 2,297 & 4,518 & 0,909 & 0,635 \\
      US & 25,635 & 5,629 & -- & 2,962 & 0,531 \\
      Other & 17,443 & 17,727 & 2,521 & 2,841 & 7,082 \\
      India & 15,275 & 0,768 & 0,148 & -- & -- \\
      Mozambique & 13,292 & -- & 2,059 & -- & -- \\
      Russia & 11,707 & 22,043 & 3,711 & 3,985 & 10,611 \\
      China & 7,034 & 3,712 & 2,149 & 2,985 & 2,664 \\
      Colombia & 2,805 & -- & -- & 0,823 & -- \\
      Singapore & 2,726 & 2,389 & 0,420 & -- & -- \\
      Philippines & -- & 6,757 & -- & -- & -- \\
      Canada & -- & 2,935 & 0,469 & 2,134 & 1,716 \\
      Vietnam & -- & -- & 0,501 & 0,509 & 0,215 \\
      Japan & -- & -- & -- & 0,294 & -- \\
      South Korea & -- & -- & -- & -- & 0,351 \\
      Peru & -- & -- & -- & -- & 0,349 \\
    \end{tabular}
  }
  \resizebox{.45\columnwidth}{!}{%
    \begin{tabular}[t]{lrrrrr}
      {\bf Exports} & Indonesia & Australia & SAFR & US & Russia \\
      India & 165,982 & 81,776 & 52,145 & 25,588 & 11,623 \\
      China & 164,523 & 31,176 & 2,283 & 5,636 & 22,045 \\
      Other & 39,718 & 29,586 & 22,092 & 28,025 & 20,836 \\
      Thailand & 26,347 & -- & -- & -- & -- \\
      Vietnam & 26,197 & 19,053 & 4,500 & -- & 3,676 \\
      Philippines & 21,150 & -- & -- & -- & -- \\
      South Korea & 15,100 & 8,410 & -- & -- & 10,617 \\
      Bangladesh & 12,736 & -- & 1,605 & -- & -- \\
      Japan & 11,334 & 19,410 & -- & 2,934 & 3,987 \\
      Taiwan & 10,499 & 11,129 & -- & -- & 3,232 \\
      Malaysia & 7,486 & 3,372 & -- & -- & -- \\
      Australia & -- & 5,819 & -- & -- & -- \\
      Brazil & -- & 6,715 & 1,214 & 11,155 & 2,175 \\
      Indonesia & -- & 6,017 & -- & -- & -- \\
      Senegal & -- & -- & 1,586 & -- & -- \\
      Pakistan & -- & -- & 6,789 & -- & -- \\
      Yemen & -- & -- & 1,234 & -- & -- \\
      Turkey & -- & -- & 2,449 & 4,249 & 8,216 \\
      UAE & -- & -- & 1,859 & -- & -- \\
      Netherlands & -- & -- & -- & 8,160 & -- \\
      US & -- & -- & -- & 2,854 & -- \\
      Morocco & -- & -- & -- & 2,140 & -- \\
      Egypt & -- & -- & -- & 5,086 & 1,613 \\
      Argentina & -- & -- & -- & 1,848 & -- \\
      Italy & -- & -- & -- & -- & 1,281 \\
    \end{tabular}
  }
  \caption{The table shows average yearly imports and exports in the dataset, for the period from January 2, 2015, to May 17, 2025. SAFR stands for South Africa.}
  \label{tab:yearly:exports:imports}
\end{table}

\begin{figure}
  \centering
  \includegraphics[width=.99\linewidth]{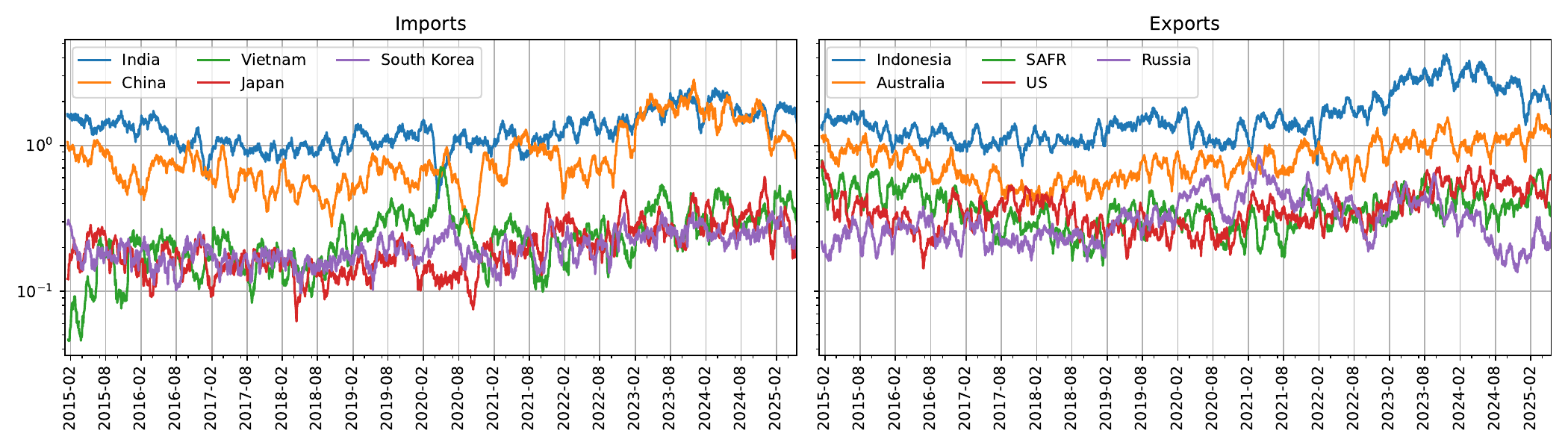}
  \caption{\emph{From left to right}: Total imports and exports to and from the five considered countries,  expressed in log scale. SAFR stands for South Africa.}
  \label{fig:time:history}
\end{figure}

\begin{figure}
  \centering
  \includegraphics[width=.47\linewidth]{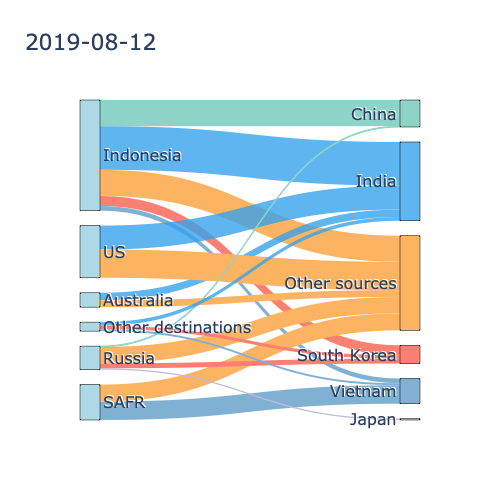}
  \includegraphics[width=.47\linewidth]{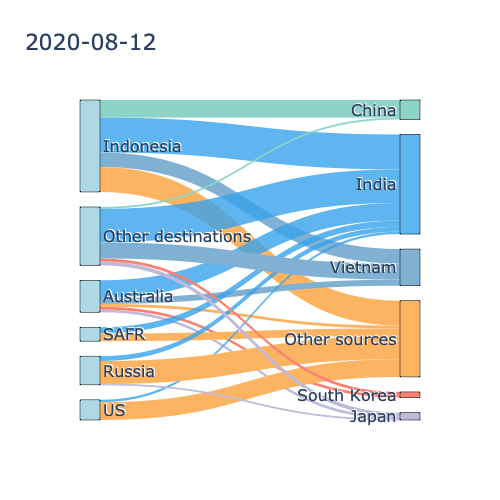}
  \caption{Sankey flow diagrams showing the impact of the pandemic on flow distributions, comparing pre-pandemic (left panel) and pandemic (right panel) periods. SAFR stands for South Africa.}
  \label{fig:sankey:pandemic}
\end{figure}

\indent We begin with an extraction of the \texttt{ShipFix} dataset, which consists of 416,977 data points detailing the imports and exports of ``Dry Coal" from January 2, 2015, to May 17, 2025, for five major exporting countries and five major importing countries, as described in Table \ref{tab:desc:exports:imports}. The latter table is reported, pictorially, in Figure \ref{fig:time:history}. The impact of the pandemic crisis (years 2020-2021) is evident in both the panels of Figure \ref{fig:time:history}, showing a decrease, particularly in the left panel representing imports, which is displayed on a logarithmic scale. The differences in trade flows before and during the pandemic are further illustrated in Figure \ref{fig:sankey:pandemic}, via two Sankey flow diagrams, which compare the structure and magnitude of flows before and during the COVID-19 pandemic. In each diagram, the width of each flow is proportional to its relative contribution to the total; this allows for a direct visual comparison of how flow magnitudes evolved between the two periods.

\begin{table}
  \centering\small\rowcolors{2}{tableShade}{white}
  \begin{tabular}{llrrr}
    &  & Number of dates & Fraction of non zero days & Number of data points \\
    Direction & Country &  &  &  \\\hline
    {\bf Exporting} & Indonesia & 3,789 & 0.48 & 41,679 \\
    & Australia & 3,789 & 0.31 & 41,679 \\
    & SAFR & 3,789 & 0.18 & 41,679 \\
    & US & 3,789 & 0.20 & 41,679 \\
    & Russia & 3,789 & 0.27 & 41,679 \\\hline
    {\bf Importing} & India & 3,795 & 0.41 & 41,745 \\
    & China & 3,794 & 0.27 & 41,734 \\
    & Vietnam & 3,791 & 0.14 & 41,701 \\
    & Japan & 3,791 & 0.14 & 41,701 \\
    & South Korea & 3,791 & 0.14 & 41,701 \\
  \end{tabular}
  \caption{Descriptive statistics for the five major exporting countries and five major importing countries, computed on a dataset of 416,977 data points. The covered period is from January 2, 2015, to May 17, 2025. SAFR stands for South Africa.}
  \label{tab:desc:exports:imports}
\end{table}

\indent Importantly, the considered countries have been strategically chosen to create a robust and comprehensive import-export matrix, despite the absence of a systematic study. The ``Dry Coal" data for these countries encompass a powerful combination of five major importing and five prominent exporting countries, ensuring a well-rounded dataset that minimizes sparsity. This dual perspective of import and export data is essential for accurately representing the ``other countries" in both sourcing and destination contexts. As a result, we develop a daily matrix of imports and exports quantified in tons, offering valuable insights into global trade dynamics. The daily averages of these flows are presented in Table \ref{tab:flows:stats}, highlighting the significant trends in our data.
\subsubsection{Transportation cost.}\label{subsubsec::selection_of_ports}
\begin{figure}
  \centering
  \includegraphics[width=\linewidth]{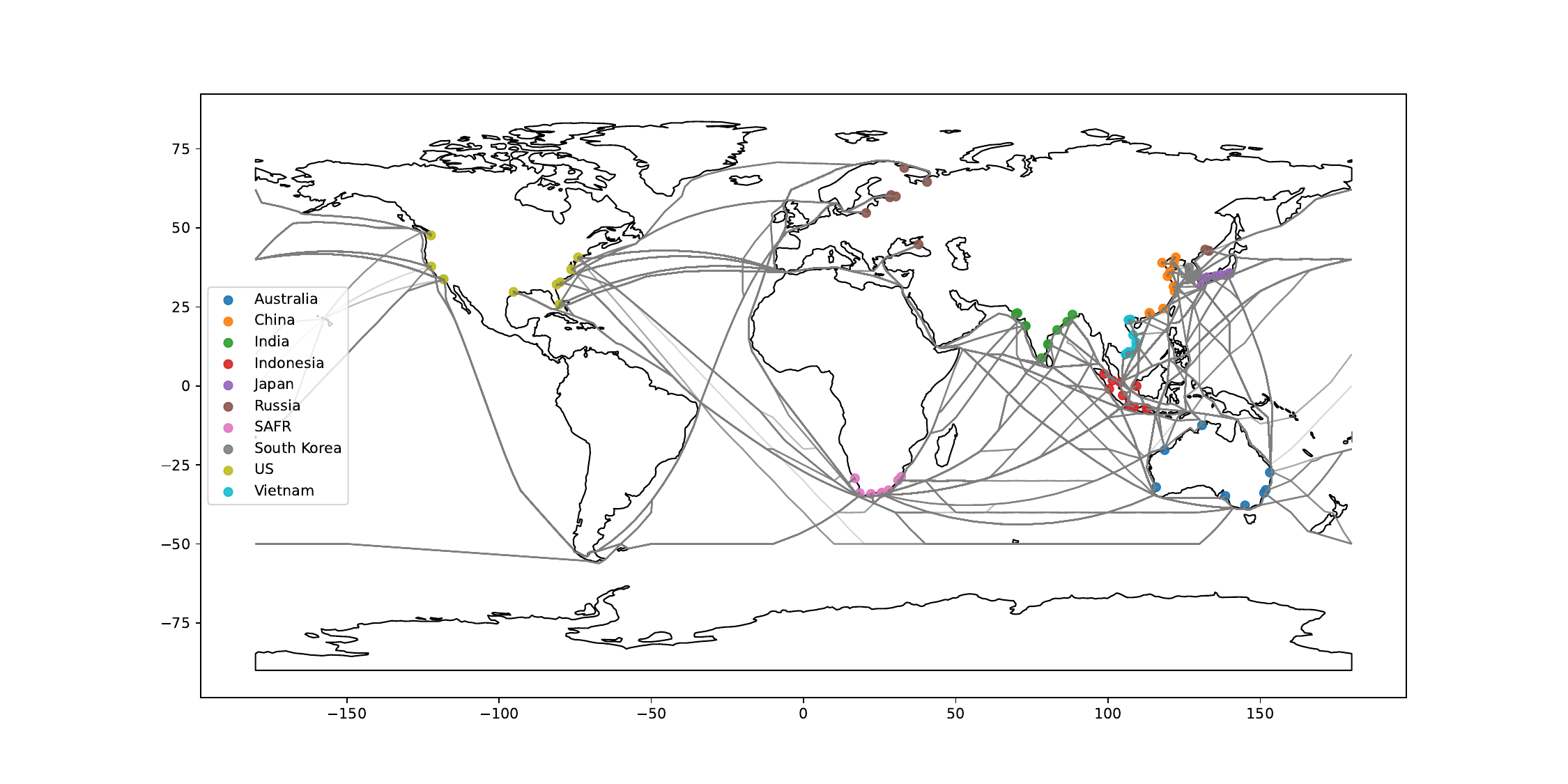}
  \caption{The figure displays the ten ports categorized by country, highlighting the greater variation in Russian and US ports compared to South African (SAFR) or Japanese ports. Additionally, all routes calculated using the \texttt{Searoutes} library are included. The shades of grey indicate the number of routes that share the same offshore segment: the darker the shade, the more routes there are.}
  \label{fig:all:routes}
\end{figure}

\begin{figure}
  \centering
  \includegraphics[width=\linewidth]{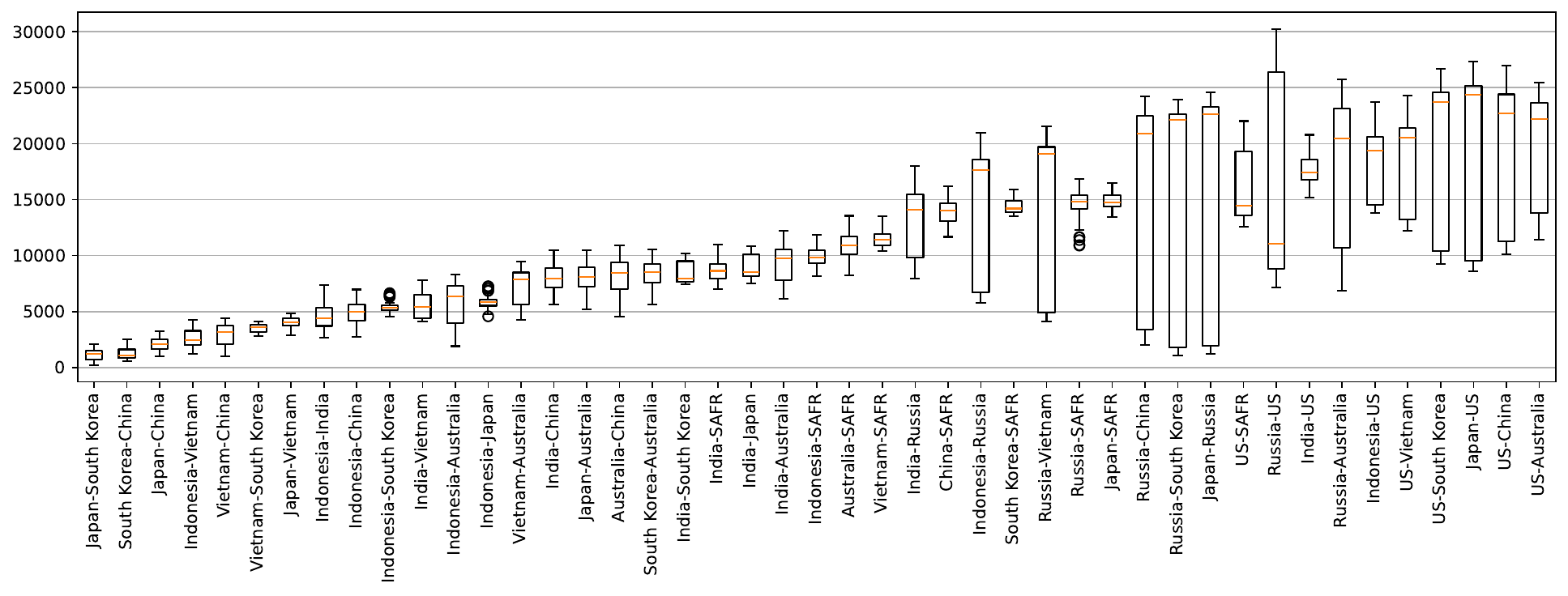}
  \caption{The dispersion of data along the main routes, organized by their averages, reveals significant differences. On the right side, we see a wide dispersion among Russian and US ports, while on the left side, South Korean and Japanese ports exhibit a much narrower concentration.}
  \label{fig:all:distances}
\end{figure}

When aggregating the import and export flows from \texttt{ShipFix} across all ports of a country, it becomes clear that the MFG model demands a transportation cost. In our model, this cost is articulated as an arbitrary function \( g \) of the transportation time \( T_{i,j} \) between the seaport \( i \) and the seaport \( j \). For this numerical application, we define \( g \) as a linear function, asserting that transportation time is directly proportional to distance. The \texttt{Searoutes} library from Eurostat provides precise maritime distances between two seaports based on their latitude and longitude (see \citep{searoute}). We conducted a thorough analysis by compiling a list of the 10 most active ports for each country considered (see Table \ref{tab:ports}). Utilizing {\sf Searoutes}, we meticulously calculated the maritime distances for all possible combinations of the $10$ ports from country \( i \) to the $10$ seaports from country \( j \), yielding $100$ distances for each pair (as depicted in Figure \ref{fig:all:routes}). The average of these distances is established as the ``distance \( T_{i,j} \)" between the two countries. While this approach serves as a robust approximation, it is important to recognize that future empirical studies can refine and enhance this methodology. Figure \ref{fig:all:distances} illustrates the dispersion of these distances through informative boxplots. It is strikingly evident from both  Figure \ref{fig:all:routes} and Figure \ref{fig:all:distances} that this approximation is less effective for a country like Russia, which has its ten most active ports distributed across vastly different seas, in stark contrast to countries such as South Korea or Japan, where port proximity yields more consistent data. This distinction underscores the criticality of considering geographical variations in future analyses.

\begin{table}
  %\centering\small
\resizebox{.25\columnwidth}{!}{
\begin{minipage}{0.45\textwidth}
    \begin{tabular}{llrr}
\toprule
Country & Port & Latitude & Longitude \\
\midrule
Australia & Sydney & -33.859 & 151.209 \\
Australia & Newcastle & -32.927 & 151.782 \\
Australia & Adelaide & -34.844 & 138.506 \\
Australia & Darwin & -12.463 & 130.846 \\
Australia & Melbourne & -37.814 & 144.963 \\
Australia & Fremantle & -32.053 & 115.748 \\
Australia & Brisbane & -27.384 & 153.118 \\
Australia & Port Hedland & -20.320 & 118.574 \\
China & Shanghai & 31.406 & 121.489 \\
China & Ningbo-Zhoushan & 29.935 & 121.811 \\
China & Shenzhen & 22.562 & 113.898 \\
China & Lianyungang & 34.731 & 119.529 \\
China & Xiamen & 24.446 & 118.065 \\
China & Dalian & 38.933 & 121.617 \\
China & Yingkou & 40.671 & 122.235 \\
China & Tianjin & 38.971 & 117.734 \\
China & Qingdao & 36.107 & 120.320 \\
China & Guangzhou & 23.101 & 113.560 \\
India & Tuticorin (V.O.C.) Port & 8.769 & 78.135 \\
India & Visakhapatnam Port & 17.690 & 83.282 \\
India & Kolkata (Syama Prasad Mookerjee) Port & 22.546 & 88.312 \\
India & Ennore (Kamarajar) Port & 13.254 & 80.340 \\
India & Chennai Port & 13.084 & 80.282 \\
India & Deendayal (Kandla) & 23.033 & 70.217 \\
India & Mundra Port & 22.741 & 69.708 \\
India & Jawaharlal Nehru Port (JNPT) & 18.963 & 72.942 \\
India & Paradip Port & 20.317 & 86.609 \\
India & Mumbai Port & 18.952 & 72.837 \\
Indonesia & Tanjung Priok & -6.104 & 106.886 \\
Indonesia & Batu Ampar & 1.141 & 104.027 \\
Indonesia & Tanjung Perak & -7.197 & 112.733 \\
Indonesia & Belawan & 3.777 & 98.686 \\
Indonesia & Teluk Bayur & -0.974 & 100.366 \\
Indonesia & Cirebon & -6.706 & 108.570 \\
Indonesia & Dumai & 1.686 & 101.449 \\
Indonesia & Gresik & -7.165 & 112.655 \\
Indonesia & Pontianak & -0.051 & 109.334 \\
Indonesia & Boom Baru & -2.999 & 104.810 \\
Japan & Osaka & 34.643 & 135.409 \\
Japan & Nagoya & 35.089 & 136.881 \\
Japan & Kobe & 34.678 & 135.187 \\
Japan & Kagoshima & 31.537 & 130.551 \\
Japan & Yokohama & 35.454 & 139.638 \\
Japan & Chiba & 35.596 & 140.089 \\
Japan & Tokyo & 35.641 & 139.764 \\
Japan & Hiroshima & 34.365 & 132.464 \\
Japan & Kitakyushu & 33.906 & 130.876 \\
Japan & Shimizu & 35.017 & 138.488 \\
% ----
\bottomrule
\end{tabular}
\end{minipage}}
\rule{10em}{0pt}\resizebox{.25\columnwidth}{!}{\begin{minipage}{0.45\textwidth}
    \begin{tabular}{llrr}
\toprule
Country & Port & Latitude & Longitude \\
\midrule
Russia & Oust-Louga & 59.673 & 28.329 \\
Russia & Saint-Pétersbourg & 59.934 & 30.306 \\
Russia & Vladivostok & 43.106 & 131.875 \\
Russia & Mourmansk & 68.971 & 33.074 \\
Russia & Arkhangelsk & 64.540 & 40.543 \\
Russia & Nakhodka & 42.783 & 132.867 \\
Russia & Novorossiysk & 44.724 & 37.768 \\
Russia & Vostochny & 42.728 & 133.099 \\
Russia & Kaliningrad & 54.698 & 20.533 \\
Russia & Primorsk & 60.367 & 28.739 \\
SAFR & Durban & -29.872 & 31.022 \\
SAFR & East London & -33.018 & 27.912 \\
SAFR & Mossel Bay & -34.181 & 22.146 \\
SAFR & Richards Bay & -28.783 & 32.038 \\
SAFR & Port Elizabeth & -33.917 & 25.602 \\
SAFR & Cape Town & -33.918 & 18.435 \\
SAFR & Port Nolloth & -29.240 & 16.861 \\
SAFR & Ngqura & -33.795 & 25.688 \\
South Korea & Yeosu & 34.754 & 127.738 \\
South Korea & Mokpo & 34.794 & 126.389 \\
South Korea & Incheon & 37.449 & 126.617 \\
South Korea & Gwangyang & 34.917 & 127.767 \\
South Korea & Busan & 35.103 & 129.040 \\
South Korea & Ulsan & 35.519 & 129.372 \\
US & Los Angeles & 33.733 & -118.276 \\
US & Long Beach & 33.754 & -118.217 \\
US & New York \& New Jersey & 40.668 & -74.041 \\
US & Savannah & 32.081 & -81.091 \\
US & Houston & 29.739 & -95.263 \\
US & Oakland & 37.795 & -122.279 \\
US & Charleston & 32.782 & -79.926 \\
US & Seattle/Tacoma & 47.602 & -122.341 \\
US & Virginia (Norfolk Portsmouth) & 36.867 & -76.311 \\
US & Miami & 25.778 & -80.183 \\
Vietnam & Cat Lai (Saïgon/HCM) & 10.768 & 106.772 \\
Vietnam & Hai Phong & 20.850 & 106.683 \\
Vietnam & Cam Pha & 21.017 & 107.300 \\
Vietnam & Da Nang & 16.083 & 108.233 \\
Vietnam & Quy Nhon & 13.767 & 109.233 \\
Vietnam & Nha Trang & 12.250 & 109.196 \\
Vietnam & Can Tho & 10.031 & 105.769 \\
Vietnam & Vung Tau & 10.353 & 107.085 \\
Vietnam & Phu My & 10.509 & 107.023 \\
Vietnam & Cai Mep–Thi Vai & 10.520 & 107.030 \\
\bottomrule
\end{tabular}
\end{minipage}}
  %}
  \caption{List of considered ports by country and their latitude and longitude.}
  \label{tab:ports}
\end{table}

\subsubsection{The notion of margins}\label{subsubsec::margins}
The maritime flow MFG model in Subsection \ref{subsec::maritme_flow_policy} is designed to capture the commercial margin generated by transporting a unit of good \( n \) from any source country \( i \) to any destination country \( j \). These margins will be derived from comprehensive data analysis. In our framework, the margin serves as a crucial indicator of the preference for directing goods to seaport \( j \) rather than seaport \( j' \). Higher expected margins create a compelling incentive for trade, but this needs to be weighed against transportation costs and concentration costs -- where an extensive congestion at ports diminishes the allure of certain destinations. We propose to define the margin \( M_{i,j} \) for shipping goods from \( i \) to \( j \) as the difference between the ``value" of the good at destination \( j \) and the source \( i \), expressed as \( M_{i,j} = v_j - v_i \). This methodology not only simplifies our model by reducing the number of parameters but also highlights a critical relationship: \( M_{i,j} \) shares the parameter \( v_i \) with \( M_{i,j'} = v_{j'} - v_i \). It is essential to understand that all these ``values" maintain their relevance and impact, remaining invariant from the model's perspective when considered within the context of an additive constant. This innovative approach ensures a more effective and precise understanding of trade dynamics, paving the way for informed decision-making in international commerce.

\subsubsection{Expected crowdedness and selection of a stationary time period.}\label{subsubsec::stationary_time_period}
In alignment with the specifications detailed in Equation \eqref{eq:quad:Q:2:1}, it is essential to choose parameters that will effectively create a proxy for the expected crowdedness of a country, which serves as a representation of the average crowdedness of its ports. To achieve this, we must identify two key parameters: a future shift, denoted as \(s\), and a moving average, denoted as \(m\). The proxy for the expected crowdedness that coordinators at seaport \(j\) will utilize on day \(t\) is defined by the following formula 

\[
\frac{1}{m}\sum_{\tau=1}^m\sum_{i=1}^{K} \Phi_{i,j}(t+s+\tau).
\]

\begin{figure}
  \centering
  \includegraphics[width=\linewidth]{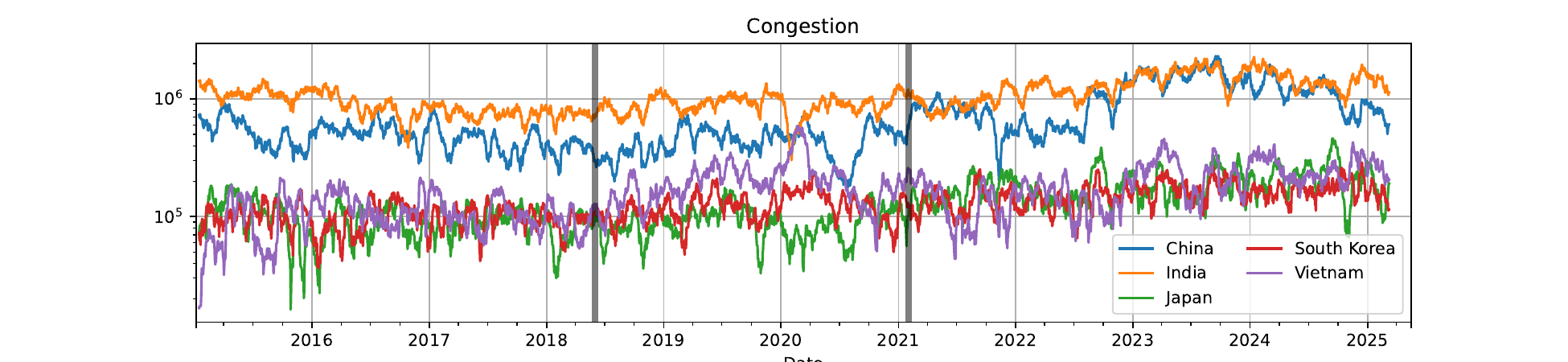}
  \caption{Time series of crowdedness for the selected countries.}
  \label{fig:crowd:ts}
\end{figure}

To establish a robust criterion, we meticulously analyzed the time series of this future crowdedness indicator and selected the parameters \(s=60\) and \(m=20\). This choice represents a 20-day moving average of the projected flows destined for the target ports over the next 60 days. Figure \ref{fig:crowd:ts} illustrates this time series data, highlighting the significance of our analysis. Notably, the period from June 2018 to January 2021 demonstrates remarkable stability, prompting us to strategically focus our inference within this timeframe in the next subsection. 

\subsection{Estimation of the parameters of the maritime flow MFG via the two-steps procedure}\label{subsec::estimation_parameter}
\begin{figure}
  \centering
  \includegraphics[width=.75\linewidth]{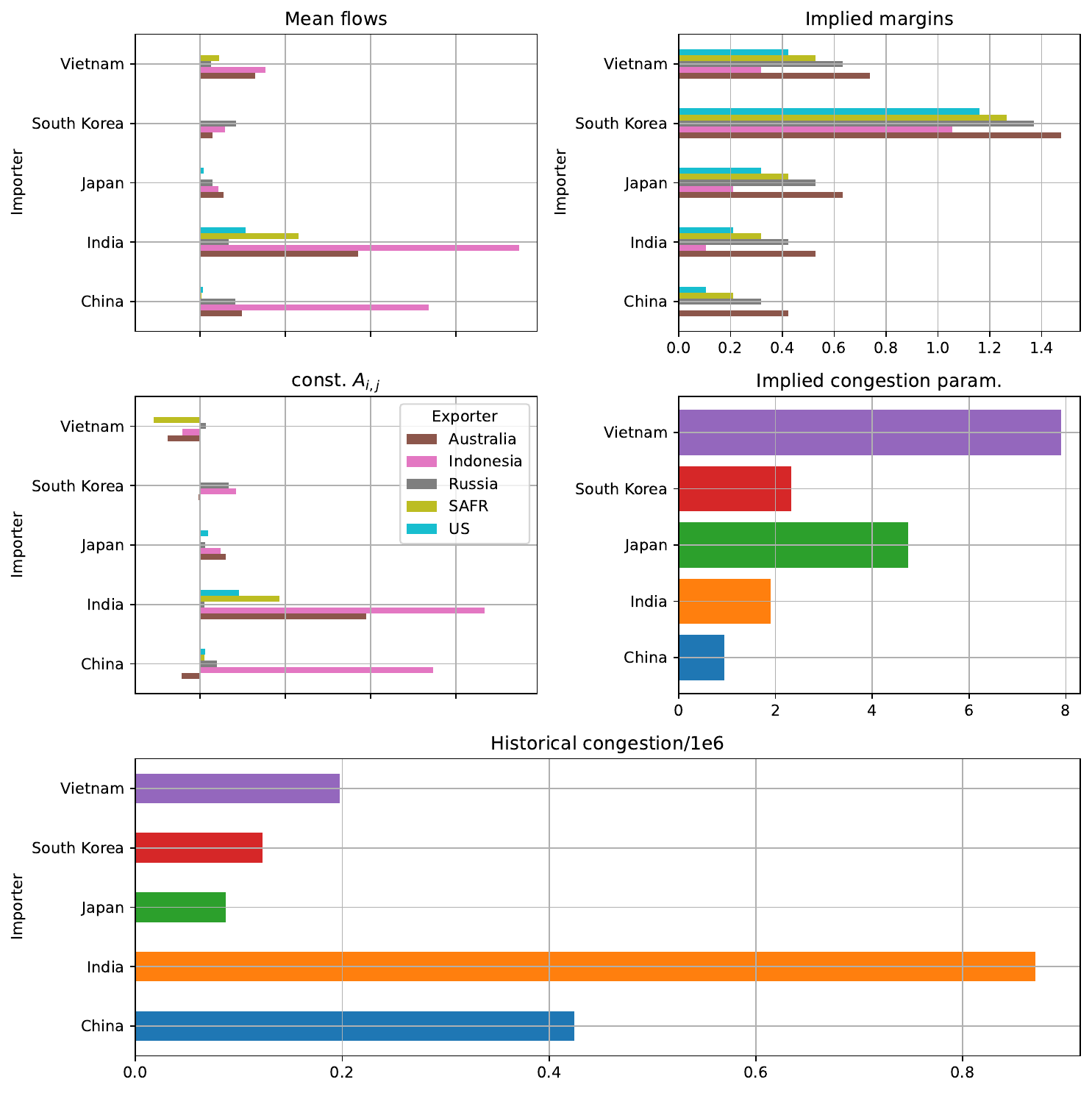}
  \caption{The figure displays the mean exports (top left), the obtained coefficients of $(A_{i,j})$ of the linear regression in Equation \eqref{eq:estimate:model} (middle left), and the mean expected congestion on destination countries $j$ (bottom panel).
  Moreover, coming from the final estimation: top right is the implied congestion coefficient $\hat r$ and middle right is the expected margin, i.e. $\hat v_j-\hat v_i$}
  \label{fig:Aij:comp}
\end{figure}

In this subsection, we validate the two-step procedure developed in Section \ref{section::statistical_inference} by demonstrating its ability to recover the parameters of the MFG model in a manner consistent with our economic understanding of coal transportation. First, for each pair of countries \((i,j)\), we estimate the coefficients \(A_{i,j}\) and $B_{i,j}$ by conducting a linear regression as outlined in Equation \eqref{eq:estimate:model}. This process involves a total of \(K(K-1)\) regressions.\\
\indent The results for the coefficients \(A_{i,j}\) are presented in Table \ref{tab:linreg:Aij}. We notice that some of the coefficients \(A_{i,j}\) are negative. This is somewhat surprising, and we now expand on this observation. First, we note that these coefficients have been defined as the expected relative margins, and are given by the following exact formula $A_{i,j} \simeq w^n_{i,j} \big(\displaystyle \wM^n_{i,j} - \sum_{\ell=1}^{K} \uw^n_{i,\ell} \,\wM^n_{i,\ell}\big).$ Then, we provide in Figure \ref{fig:Aij:comp} the average exports (\emph{top panel}), the estimated coefficients $( A_{i,j} )$ (\emph{mid panel}), and the average expected congestion in destination countries \( j \) (\emph{bottom panel}), for each pair of countries \((i,j)\). We note that the coefficients \(A_{i,j}\) correlate with the average export levels by indicating larger expected margins. However, this correlation holds true only when the average crowding coefficient is not excessively high. In particular, upon analysing Figure \ref{fig:Aij:comp}, it becomes evident that for instance Vietnam's limited ``Dry Coal" imports, combined with relatively low congestion levels compared to other nations, lead to a compelling conclusion. The model demonstrates that the expected margin must be set lower than that of other countries. Consequently, this results in a negative value for $A_{i,j}$. This insight underscores the critical relationship between trade dynamics and market conditions in shaping outcomes. 

\begin{table}
  \centering\rowcolors{2}{tableShade}{white}
  \begin{tabular}{lrrrrr}
    \hfill to & China & India & Japan & South Korea & Vietnam \\
    from &  &  &  &  &  \\
    Australia & -21,331 & 195,306 & 30,240 & -1,539 & -38,075 \\
    Indonesia & 273,369 & 333,976 & 24,056 & 42,270 & -20,895 \\
    Russia & 19,676 & 5,617 & 5,787 & 33,787 & 7,185 \\
    SAFR & 5,155 & 93,401 & -- & -- & -54,015 \\
    US & 6,510 & 45,662 & 9,589 & -- & -- \\
  \end{tabular}
  \caption{The table displays the constants $A_{i,j}$ obtained via the linear regression in Equation \eqref{eq:estimate:model}.}
  \label{tab:linreg:Aij}
\end{table}

\begin{figure}
  \centering
  \includegraphics[width=\linewidth]{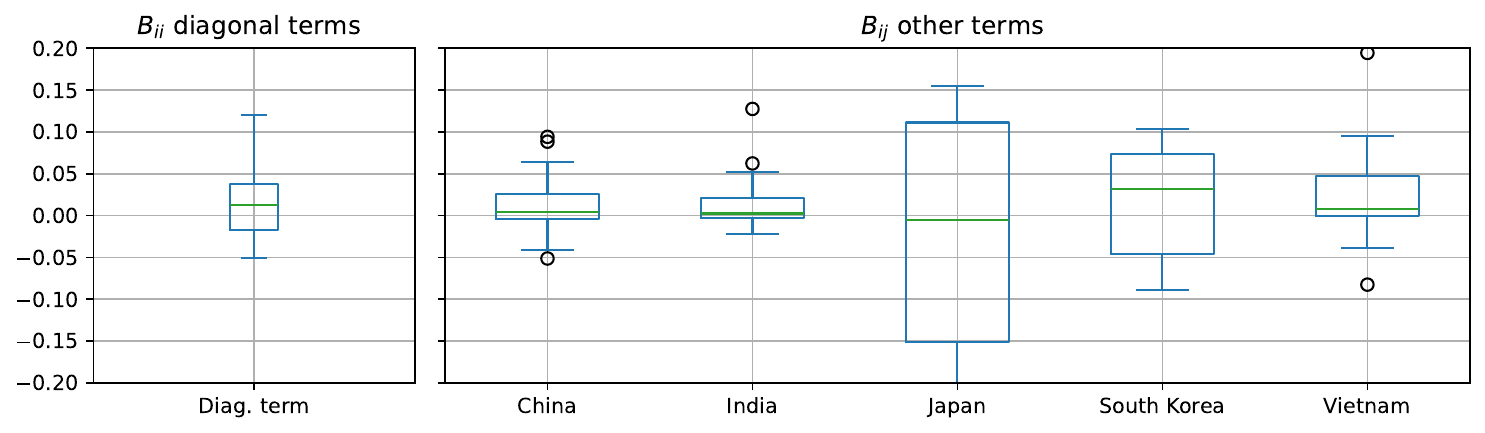}
  \caption{Distribution of the regression slope coefficients $(B_{i,j})$. The diagonal terms are reported in the left panel, whereas the non-diagonal terms are reported in the right panel.}
  \label{fig:distrib:Bij}
\end{figure}

\indent On the other hand, the regression slope coefficients $( B_{i,j} )$ are reported in the last column of Table \ref{tab:corr:B:cong}. In Figure \ref{fig:distrib:Bij}, we, instead show the distribution of the coefficients $B_{i,j}$. In the left panel, it displays the diagonal terms, whereas in the right panel it shows the non-diagonal terms, strategically grouped by destination country. Before commenting the results, it is important to highlight that a distinct linear regression analysis is conducted for each unique pair $(i,j)$, emphasizing the robustness and specificity of the findings.\\
\indent The theoretical formula for the slope coefficients is as follows:
\begin{equation*}
    B_{i,j} = w^n_{i,j} \sum_{\ell=1}^{K} (\uw^n_{i,\ell} - \one_{\{\ell=j\}}) \, r^n_\ell.
\end{equation*}
In particular, we do not anticipate a perfect equality due to the inherent non-stationarity of the real data. Our analysis is constrained to a three years time period (see Subsection \ref{subsubsec::stationary_time_period}), 
yielding only 958 data points for each linear regression model. Despite this limitation, we find that all the $(\uw^{n}_{i}$ values are less than one, leading us to expect that most coefficients will be positive -- except for the diagonal term (i.e., \(B_{i,i} \leq 0\) and \(B_{i,j} \geq 0\) for \(i \neq j\)). Significantly, the diagonal terms are notably smaller than their non-diagonal counterparts, as illustrated by the contrasting $y-$values in the data presented in the right and left panels of Figure \ref{fig:distrib:Bij}.  While the medians of the non-diagonal coefficients in our data sample tend to be positive, we must acknowledge that not all are; particularly for Japan, this finding raises red flags. It implies that the volume of ``Dry Coal" exported to Japan is generally not positively correlated with the expected congestion there; in fact, the correlation may be negative. Table \ref{tab:corr:B:cong}, third column, substantiates this claim, showing that all negative non-diagonal terms correspond to such adverse relationships, with the sole exception of the route from Russia to North Korea.  This evidence suggests that the constants we have selected -- namely, a 60-day lead time and an average over 20 days -- might be inadequate for accurately forecasting transport times to Japan. Clearly, this warrants further investigation to refine our approach and enhance the accuracy of our predictions.\\
\indent Turning now to the second step of our statistical methodology, we run the following minimization problem   
\begin{equation*}\label{eq::minimization}
   \min_{\substack{c_n, (r_j,v_j)_{1 \leq j \leq K}}}\sum_{\substack{i,j=1 \\ i \neq j}}^{K} \left\|A_{i,j} \cdot (r_j+c_ng(T_{i,j})\one_{\{\ell=j\}}) - (v_j-v_i) \right\|^{2}.
\end{equation*}
In particular, this means that we effectively integrate the results from the $K(K-1)$ linear regressions, utilizing their coefficients \(A_{i,j}\) to create a cohesive analysis. Our connections are grounded in the following criteria:
\begin{itemize}
    \item We establish consistent source and destination terminal values \(v_i\) and \(v_j\). For instance, the cost of sending coal to Japan remains uniform whether it originates from Australia or the US, ensuring that both linear regressions operate with the same foundational values.
    \item We assert that the congestion cost coefficient \(c_j\) at any given destination is invariant, providing clarity and consistency across different source countries.
    \item  We apply a uniform coefficient \(c\) to the distance between each source \(i\) and destination \(j\). This approach guarantees that all variability is thoroughly accounted for by the distance factor alone, enhancing the robustness of our analysis.
\end{itemize}

\begin{figure}
  \centering
  \includegraphics[width=.8\linewidth]{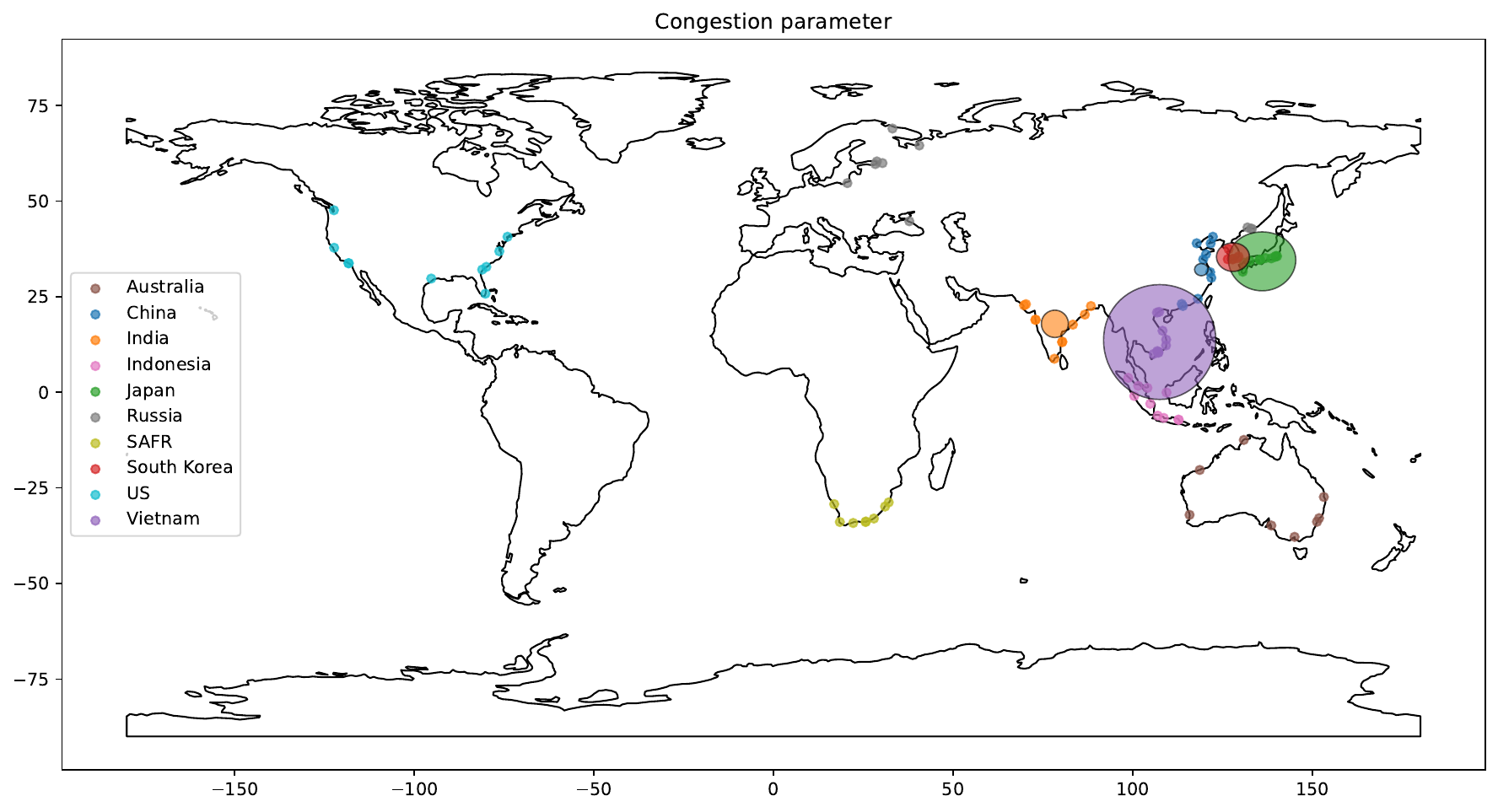}
  \includegraphics[width=.8\linewidth]{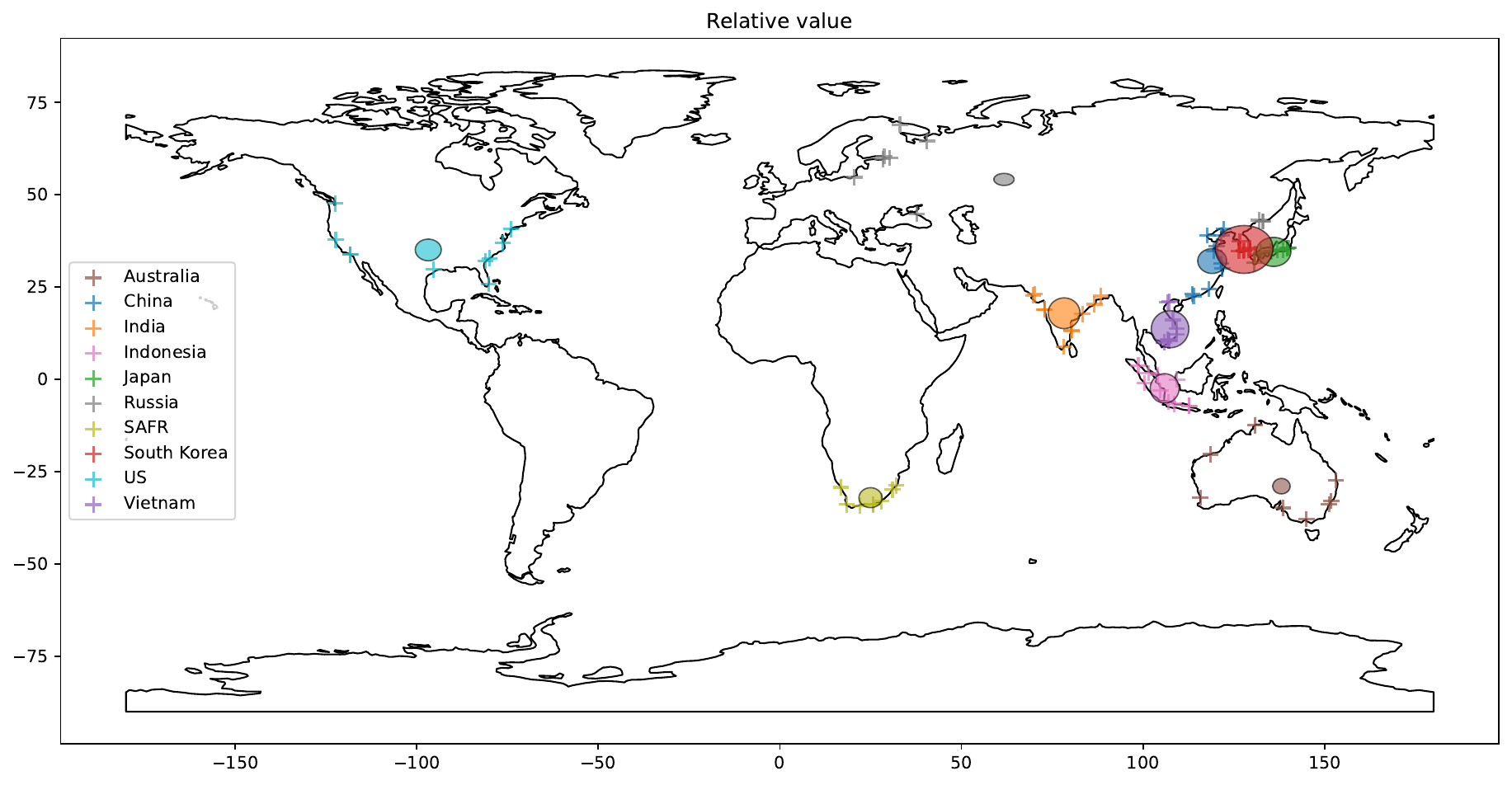}
  \caption{The figure shows a world map of the congestion coefficients (top panel) and the expected relative value of ``Dry Coal" (bottom panel) estimated by the MFG model using the \texttt{ShipFix} dataset at the country level. SAFR stands for South Africa.}
  \label{fig:worldmap}
\end{figure}

The results can be found in Table \ref{tab:result:crv} and Figure \ref{fig:worldmap}. 
To read the table (or the world maps, that contains the same information), one can start with the lower panel of  \ref{fig:worldmap}:
The diameter of the circles are proportional to the estimate of the value of the good (here ``Dry Coal") at the country.
It means agents of the mean field (ships) have a direct interest to make margin to transport from smaller circles to largest ones:
Figure \ref{fig:Aij:comp} is useful to comment these implied margins since its top right bar chart is made of the estimated $v_j-v_i$.
The reader can see that combining this intensive to ship goes more to Vietnam, South Korea and Japan than the historical imports (top left panel of the same figure), but that the estimated congestion coefficients (middle right panel) goes in the other direction, in particular for Vietnam and Japan. It ``explains'' why flows are not going that much to these countries but rather to India and China.
It is noticeable that the implied congestion parameter (middle right) is not aligned with the historical time series of congestion (the average values are in the bottom panel): the model takes the historical congestion flows into account and the congestion parameter is only a multiplier of the congestion. The effective congestion cost for the model is the parameters $r$ multiplied by the squared expected traffic (i.e. congestion) in the given country.

The reader should keep in mind this numerical application is at the level of the countries (using the sum of the traffic in the top ports of each country) and it is expected that some accuracy is lost during such a macroscopic approximation. 
A fine atomic study is possible but goes beyond the scope of this section of this paper that is there for illustration, to demonstrate that the proposed approximation methodology is numerically realistic.

A estimations of the parameters of the MFG thanks to our results and using our proposed methodology deserves a full research by itself, at the level of ports, on multiple goods, and should target to go beyond the non stationarity of the data by capturing slow varying preferences of the agents. 
It would open the door to ``What if scenarios'' particularly suited to understand the consequences of physical or geopolitical shocks, like wildfires, wars, or sanctions.

\begin{table}
  \centering\rowcolors{2}{tableShade}{white}
\begin{tabular}{llrr}
country from & country to & corr & $B_{i,J}$ \\
Indonesia & Japan & -0.008 & 0.016 \\
Indonesia & Vietnam & 0.095 & 0.042 \\
Indonesia & South Korea & -0.043 & -0.038 \\
Indonesia & China & -0.059 & 0.001 \\
Indonesia & India & -0.020 & -0.045 \\
Australia & Japan & -0.016 & 0.121 \\
Australia & Vietnam & 0.106 & 0.075 \\
Australia & South Korea & 0.072 & 0.031 \\
Australia & China & 0.169 & 0.101 \\
Australia & India & -0.050 & -0.044 \\
SAFR & Vietnam & 0.178 & 0.097 \\
SAFR & India & -0.004 & 0.031 \\
SAFR & China & -0.006 & -0.003 \\
US & Japan & -0.038 & -0.013 \\
US & India & -0.033 & -0.016 \\
US & China & -0.012 & -0.017 \\
Russia & Japan & -0.016 & -0.047 \\
Russia & Vietnam & 0.062 & 0.012 \\
Russia & South Korea & 0.033 & -0.051 \\
Russia & China & 0.096 & 0.038 \\
Russia & India & 0.027 & 0.013 \\
\end{tabular}
  \caption{The table displays the regression slope coefficients $(B_{i,j})$ and the correlation between the volume of ``Dry Coal" exported to a country with the expected congestion on that country.}
  \label{tab:corr:B:cong}
\end{table}

\begin{table}
  \centering\rowcolors{2}{tableShade}{white}
  \begin{tabular}{lrrr}
    & cong. costs $r$ & $v_j$ & $v_i$ \\\hline
China & 1.05 & 1.16 & -- \\
India & 2.10 & 1.28 & -- \\
Japan & 5.24 & 1.40 & -- \\
South Korea & 2.56 & 2.33 & -- \\
Vietnam & 8.73 & 1.51 & -- \\
Australia & -- & -- & 0.70 \\
Indonesia & -- & -- & 1.16 \\
Russia & -- & -- & 0.81 \\
SAFR & -- & -- & 0.93 \\
US & -- & -- & 1.05 
\end{tabular}
  \caption{The table displays the results of minimizing the transportation and congestion costs, as well as the expected ``Dry Coal" value at both the source and destination, derived from the coefficients of the $K(K-1)$ linear regressions. The calculated transportation costs amounted to $c=0.92.$}
  \label{tab:result:crv}
\end{table}

\begin{table} % TODO: reorder 'others'
  \centering\small\rowcolors{2}{tableShade}{white}
 \begin{tabular}{llrrrrrrr} % {1.05\textwidth}
 & country\_to & China & India & Japan & South Korea & Vietnam & Others & All \\\hline
country\_from & stat &  &  &  &  &  &  &  \\
\multirow[c]{1}{*}{Australia} & mean & 85,437 & 224,137 & 53,220 & 23,037 & 52,240 & 171,710 & 609,780 \\
 & std & 155,595 & 247,722 & 100,049 & 61,589 & 98,260 & 199,166 & -- \\
\multirow[c]{1}{*}{Indonesia} & mean & 451,091 & 455,033 & 31,076 & 41,402 & 71,829 & 323,257 & 1,373,687 \\
 & std & 556,317 & 460,926 & 61,568 & 78,164 & 115,554 & 322,486 & -- \\
\multirow[c]{1}{*}{Russia} & mean & 60,446 & 31,870 & 10,932 & 29,110 & 10,079 & 102,401 & 244,837 \\
 & std & 105,251 & 62,377 & 30,881 & 52,818 & 30,999 & 128,298 & -- \\
\multirow[c]{1}{*}{SAFR} & mean & 6,259 & 142,933 & -- & -- & 12,338 & 106,399 & 267,930 \\
 & std & 28,362 & 186,638 & -- & -- & 45,604 & 128,452 & -- \\
\multirow[c]{1}{*}{US} & mean & 15,454 & 70,159 & 8,045 & -- & -- & 174,036 & 267,693 \\
 & std & 49,518 & 119,521 & 31,008 & -- & -- & 196,374 & -- \\
\multirow[c]{1}{*}{Others} & mean & 93,230 & 161,339 & 28,681 & 37,036 & 23,502 & -- & 343,788 \\
 & std & 116,169 & 178,432 & 56,806 & 63,522 & 49,103 & -- & -- \\
\multirow[c]{1}{*}{Total} & mean & 711,917 & 1,085,471 & 131,953 & 130,584 & 169,988 & 877,802 & 3,107,716 \\
 %& std & -- & -- & -- & -- & -- & -- & -- \\
\end{tabular}
  \caption{The table displays descriptive statistics -- mean and standard deviation (std) -- on the average daily flows, as measured in tons. SAFR stands for South Africa.}
  \label{tab:flows:stats}
\end{table}

\newpage
\section{Conclusions and future research}\label{section::conclusion_and_future_research}
The goal of this article was to advance the MFG framework in seaport competition, which is still a nascent field that entails a lot of open questions and challenges. To this end, we have developed a MFG model to understand maritime flows and the sources of decision-making. We have considered transportation costs, expected profit margins from loading specific goods at seaports and unloading them at various destinations, and a congestion factor that represents the cost of accessing the destination port.  We have provided a clear solution to the stationary version of the MFG model under the assumption of quadratic costs and that a sufficient number of ships are involved, so that the stochastic transition matrix between seaports closely approximates its deterministic equivalent. In particular, we have found that for each category of goods, the MFE distribution of trade flows is governed by the first eigenvector associated with the system's optimal control variables. These eigenvectors emerge as solutions to a set of interdependent linear equations that share a common mean field, capturing the aggregate influence of all interacting flows within the network. Furthermore, we have established a unifying condition that integrates all relevant cost parameters --- transportation, handling, and transition costs --- to guarantee the existence and uniqueness of this equilibrium solution. We have developed a statistical methodology to infer the parameters of the game, i.e., the costs and the commercial margins’ components, from the \texttt{ShipFix} dataset, which does not require the application of, e.g., fictitious play techniques to stabilize its solution. As good, we have taken into account the ``Dry Coal". Finally, We have analyzed the results to clarify how the various cost components are influenced by the empirical flows.\\
\indent Much remains to be done. From a theoretical perspective, the most restrictive assumption is that we equate the stochastic realization of a Markov Chain -- representing the desired transitions of ships between ports -- with its average. This assumption holds reasonably well when there is a very large number of ships; however, future theoretical work could relax this condition. Regarding the statistical methodology, although our results appear ``validated," or at least compatible, with the established economic knowledge regarding ``Dry Coal" transport, there are a few assumptions that could be reconsidered in future research. These include conducting a similar estimation at the level of individual ports and incorporating explicit production terms.\\
\indent Nonetheless, we believe this work represents a significant step toward a more refined application of MFGs in maritime traffic management.

% Acknowledgments here
\ACKNOWLEDGMENT{Authors would like to thank ShipFix for the data; Evan Show for the early empirical work on maritime traffic; Michele Bergami,  Simone Moawad, and Barath Raaj Suria Narayanan for the explorations they have made during their Master Capstone project at The London School of Economics and Political Science; and Borel Martial Domgue Defo for the work on the ShipFix dataset on Coal during his internship at the CMAP, Ecole Polytechnique.}% Leave this (end of acknowledgment)

% Appendix here
% Options are (1) APPENDIX (with or without general title) or 
%             (2) APPENDICES (if it has more than one unrelated sections)
% Outcomment the appropriate case if necessary
%
% \begin{APPENDIX}{<Title of the Appendix>}
% \end{APPENDIX}
%
%   or 
%
% \begin{APPENDICES}
% \section{<Title of Section A>}
% \section{<Title of Section B>}
% etc
% \end{APPENDICES}

% References here (outcomment the appropriate case) 
\newpage
% CASE 1: BiBTeX used to constantly update the references 
%   (while the paper is being written).
\bibliographystyle{informs2014trsc} % outcomment this and next line in Case 1
\bibliography{mfg_ships} % if more than one, comma separated

% CASE 2: BiBTeX used to generate mypaper.bbl (to be further fine tuned)
%\input{mypaper.bbl} % outcomment this line in Case 2

\end{document}